\documentstyle[12pt,amssymb]{amsart}
\newenvironment{Defn}{\par\smallskip\noindent{\it Definition.\/}}
{\par\medskip\endtrivlist}
\newenvironment{Thm}{\par\smallskip\noindent{\it Theorem.\/}}
{\par\medskip\endtrivlist}
\newenvironment{Rem}{\par\smallskip\noindent{\it Remark.\/}}
{\par\medskip\endtrivlist}

\numberwithin{equation}{section}
\theoremstyle{plain} 
\newtheorem{thm}{Theorem}[section]
\newtheorem{lem}[thm]{Lemma}
\newtheorem{prop}[thm]{Proposition}
\newtheorem{cor}[thm]{Corollary}
\theoremstyle{definition}
\newtheorem{exmp}[thm]{Example}
\newtheorem{defn}[thm]{Definition}  

\def\H{\bold{H}}   
\def\K{\bold{K}}   
\def\Lat{{\cal{L}}}
\def\N{{\cal{N}}}
\def\~{\tilde}
\def\Omega{\varOmega}
\def\<{\langle}
\def\>{\rangle}
\def\AB{A\text{-}B}
\def\oEnd{\operatorname{End}}
\def\Aut{\operatorname{Aut}}

\begin{document}

\bibliographystyle{amsplain}

\title{Normal intermediate subfactors}
\author{Tamotsu Teruya}
\address{Department of Mathematics \\
  Hokkaido University \\
  Sapporo 060 Japan}
\email{t-teruya@@math.hokudai.ac.jp}
\keywords{subfactor, Jones index, normal, Hopf algebra, Kac
  algebra,
  strongly outer action,  modular lattice}
\date{Oct.~30, 1995}
\maketitle
\begin{abstract}
  Let $N \subset M$ be an irreducible inclusion of type 
  type II$_1$ factors with finite Jones index.
  We shall introduce the notion  of normality for intermediate 
  subfactors of the inclusion $N \subset M$. 
  If the depth of $N \subset M$ is 2, then 
  an intermediate subfactor $K$ for $N \subset M$ is normal in 
  $ N \subset M$ if and only if the depths of $N \subset K$ 
  and $K \subset M$ are both 2.
  In particular, 
  if $M$ is the crossed product $N \rtimes G$ of a finite group
  $G$, 
  then $K = N \rtimes H$ is normal in $N \subset M$ if and 
  only if $H$ is a normal subgroup of $G$.
\end{abstract}

\section{Introduction}
\label{introduction}

The index theory for type II$_1$ factors initiated by 
V.~Jones \cite{Jo:index} and
the classification of subfactors has been studied by 
many people (\cite{EK:orbifold}, \cite{Izumi:fusion},
\cite{Izumi:index3}, \cite{IK:Dn},  
\cite{Kawahi:flat}, \cite{Ko:charac}, \cite{KL:remini}, 
\cite{Loi:typeIII}, \cite{Lo:index,Lo:indexII},
\cite{PP:entropy}, ...). 
A.~Ocneanu 
\cite{Oc:Qg} introduced the concept of 
a paragroup 
to classify subfactors.
By using the so-called standard invariant equivalent to 
the paragroups, 
S.~Popa \cite{Po:amenable},\cite{Po:reduction}
classified subfactors under more general conditions.
Paragroup or the standard invariant 
for an inclusion of type II$_1$ 
factors with finite Jones index is a group like object 
which contains finite groups. 
So the theory of finite groups may be considered as
part of the index theory for 
an inclusion of type II$_1$ factors with finite 
Jones index.
It is  well known that if $\alpha :G \to Aut(N)$ is 
an outer action of a finite group $G$ on a type II$_1$ factor 
$N$ and $K$ is an intermediate subfactor for 
$N \subset N \rtimes_{\alpha}G$, then there is a subgroup 
$H$ of $G$ such that $K = N \rtimes_{\alpha} H$ 
(see for instance \cite{NT:Galois}). 
On the other hand, Y.~Watatani \cite{Wa:lattice} showed that 
there exist only finitely many intermediate subfactors 
for an irreducible inclusion with finite index. 
So it is natural to consider 
intermediate subfactors    
as  ``quantized subgroups'' in the index theory 
for an inclusion of type II$_1$ factors. 
The notion of normality for subgroups plays important role 
in the theory of 
finite groups. 
In this note we introduce the notion of 
normality for intermediate subfactors of 
irreducible inclusions.

D.~Bisch \cite{Bi:interm} and 
A.~Ocneanu \cite{Oc:QS} gave a nice 
characterization of intermediate 
subfactors of a given 
irreducible inclusion $N \subset M$ 
in terms of Jones projections and  
Ocneanu's Fourier transform 
$\cal F: N' \cap M_1 \to M' \cap M_2$.
We define normal intermediate subfactors as follows:
\begin{Defn}
  Let $N \subset M$ be an irreducible inclusion of type II$_1$ 
  factors with finite index and $K$ an intermediate subfactor 
  of the inclusion $N \subset M$.
  Then $K$ is a {\it normal intermediate subfactor\/} of 
  the inclusion $N \subset M$ if 
  $e_K \in \cal Z(N' \cap M_1)$ and 
  $\cal F(e_K) \in \cal Z(M'\cap M_2)$, 
  where $e_K$ is the Jones projection for the inclusion 
  $K \subset M$.
\end{Defn}

Every finite dimensional Hopf $C^*$-algebra (Kac algebra)
gives rise to an  irreducible inclusion of AFD II$_1$ factors, 
which are characterized by depth 2 
(see for example \cite{Oc:QS}, \cite{Sek:Kac}, \cite{Szym:Hopf}, 
\cite{Yuchi:Kac}).
Let $M$ be the crossed product algebra $N \rtimes \H$ of 
$N$ by an outer action of a finite 
dimensional Hopf $C^*$-algebra $\H$.
Unfortunately, there is no one-to-one correspondence between 
the intermediate subfactors of $N \subset M$ and 
the subHopf $C^*$-algebras of $\H$ in general.
But we get the next result:
\begin{Thm}
  Let $N \subset M$ be an irreducible, depth 2 inclusion 
  of type II$_1$ factors with finite index, i.e., 
  $M$ is described as the crossed product algebra 
  $N \rtimes \H$ of $N$ by an outer action of a finite 
  dimensional Hopf $C^*$-algebra $\H$.
  Let $K$ be an intermediate subfactor of $N \subset M$ 
  and $e_K$ is the Jones projection for $K \subset M$.
  Then $K$ is described as the crossed product algebra 
  $N \rtimes \K$ of $N$ by an outer action of a
  subHopf $C^*$ algebra $\K$ of $\H$ 
  if and only if 
  $e_K$ is an element of the center of the relative 
  commutant algebra $N' \cap M_1$, 
  where $M_1$ is the basic extension for $N \subset M$.
\end{Thm}
Let $N \subset M$ be an irreducible inclusion of 
type II$_1$ factors with finite index and 
$M_1$ the basic extension 
for $N \subset M$.
Let $K$ be an intermediate subfactor  of $N \subset M$ 
and $K_1$ the basic extension for $K \subset M$.
Then $K_1$ is an intermediate subfactor of $M \subset M_1$.
For the Jones projections $e_K$ and $e_{K_1}$ for 
the inclusions 
$K \subset M$ and $K_1 \subset M_1$, respectively,
since $\cal F(e_K) = \lambda e_{K_1}$ for some scalar 
$\lambda$, 
we get the next theorem: 
\begin{Thm}
If the depth of a given irreducible 
inclusion $N \subset M$ is 2, then 
an intermediate subfactor $K$ of $N \subset M$ is 
normal in $N \subset M$ if and only if 
the depths of $N \subset K$ and 
$K \subset M$ are both 2. 
\end{Thm}
The author \cite{Tamo:normal} showed that 
if $M$ is the crossed product $N \rtimes G$ of 
finite group $G$ and $K = N \rtimes H$, then 
$H$ is a normal subgroup of $G$ if and only  
if $K \subset M \simeq K \subset K \rtimes F$ 
for some finite group $F$, i.e., the depth of 
$K \subset M$ is 2.
Hence we have $H$ is a normal subgroup of $G$ if and only if 
$K$ is a normal intermediate subfactor of $N \subset M$
by the previous theorem.
Therefore our notion of normality for intermediate subfactors 
is an extension of that in the theory of 
finite groups.

\bigskip
\noindent
{\bf Acknowledgments}.
I should like to express my gratitude to Professor 
H.~Kosaki for helpful advice and suggestion about examples of 
normal intermediate subfactors and non normal ones. 
(group type inclusion, strongly outer action, ...).
And I should like to  thank
Professor Y.~Watatani for fruitful discussions,
may useful comments  
and constant encouragement.
I also thank Professor Y.~Sekine and Professor T.~Yamanouchi
for useful comments.

\section{preliminaries}
\label{sec:pre}

\subsection{intermediate subfactors}
\label{sec:inter}

We recall here some results for intermediate subfactors. 
Let $N \subset M$ be a pair of type II$_1$ factors.
We denote by $\Lat(N \subset M)$ the set of all 
intermediate von Neumann subalgebras of $N \subset M$.
The set $\Lat(N \subset M)$ forms a lattice under the two 
operations $\vee$ and $\wedge$ defined by 
\begin{displaymath}
   K_1 \vee K_2 = (K_1 \cup K_2)'' \ \text{and} \ 
   K_1 \wedge K_2 = K_1 \cap K_2.
\end{displaymath}
If the relative commutant algebra $N' \cap M$ is trivial,
then $\Lat(N \subset M)$ is exactly the lattice of intermediate
subfactors for $N \subset M$. In fact for any 
$K \in \Lat(N \subset M)$, 
$\cal Z(K) = K' \cap K \subset N' \cap M = \Bbb C$.
If $M$ is the  crossed product $N \rtimes_{\alpha} G$
for an outer action $\alpha$ of a finite group $G$, 
then it is well known that the intermediate subfactor lattice 
$\Lat(N \subset M)$ is isomorphic to 
the subgroup lattice $\cal{L}(G)$ 
(see  \cite{NT2:Galois2}, \cite{NT:Galois}). 
In \cite{Wa:lattice} Y.~Watatani proved the next theorem.
\begin{Thm}
  Let $N \subset M$ be a pair of type II$_1$ factors. 
  If 
  $[M:N] < \infty$ and $ N' \cap M = \Bbb C$, then 
  $\Lat(N \subset M)$ is a finite lattice.
\end{Thm}
\noindent
Later we were noted  that this theorem 
was shown 
by S.~Popa implicitly \cite{Po:corres}.

From now on we assume that $[M: N] < \infty $ and 
$N' \cap M = \Bbb C$.
Let 
\begin{displaymath}
  N \subset M \subset M_1 \subset M_2
\end{displaymath}
be the Jones tower of $N \subset M$ 
obtained by iterating the basic extension.
Let 
$e_N \in M_1$ and $e_M \in M_1$ be the Jones projections 
for $N \subset M$ and $M \subset M_1$, respectively.
We denote by $\cal F$, Ocneanu's Fourier transform from 
$N' \cap M_1$ onto $M' \cap M_2$ 
i.e., 
\begin{displaymath}
  \cal F (x) 
  = [M:N]^{- \frac{3}{2}}E^{N'}_{M'}(xe_Ne_M), 
  \  x \in N' \cap M_1, 
\end{displaymath}
where $E^{N'}_{M'}$ is the conditional expectation 
from $N'$ onto $M'$.
For $K \in \Lat(N \subset M)$, if $e_K$ is the Jones projection 
for $K \subset M$, then $e_K$ is an element of $N' \cap M_1$.
In fact $K_1 = \langle M, e_K \rangle = J_M K'J_M \subset 
J_M N' J_M = M_1$ and hence $e_K \in K' \cap K_1 \subset N' \cap M_1$.

D.~Bisch \cite{Bi:interm} and A.~Ocneanu \cite{Oc:QS} gave the next
characterization of intermediate  subfactors  in terms of  
Jones projections in $N' \cap M_1$.
\begin{Thm}
  Let $p$ be a projection in $N' \cap M_1$.
  There exists an intermediate subfactor $K \in \Lat(N \subset M)$ such 
  that $p = e_K$ if and only if 
  \begin{enumerate}
  \item \begin{math}
      p \geq e_N,
    \end{math}
  \item \begin{math}
      \cal F(p) = \lambda q  \ \text{for some} \ \lambda \in \Bbb C
      \ \text{and some projection} \ q \in M' \cap M_2.  
    \end{math}
  \end{enumerate}
  In this case, $q$ is  the Jones projection $e_{K_1}$ for 
  $K_1 \subset M_1$.
\end{Thm}
For the convenience, we prove the next lemmas 
(see for example \cite{Bi:interm}, \cite{SW:angle}).
\begin{lem}\label{lem:projection}
  With the above notations, 
  we have
  \begin{displaymath}
    e_K = [K:N][M:N] E^{M_2}_{M_1}(e_Me_Ne_{K_1}),
  \end{displaymath}
  where $E^{M_2}_{M_1}$ is the trace preserving conditional
  expectation form $M_2$ onto $M_1$.
\end{lem}
\begin{pf}
  Since $e_M \leq e_{K_1}$, we have 
  \begin{displaymath}
    e_Me_Ne_{K_1} = e_Me_{K_1}e_Ne_{K_1} 
    = e_ME^{M_1}_{K_1}(e_N).  
  \end{displaymath}
  Since
  $
    E^{M_1}_{K_1}(e_N)e_K = E^{M_1}_{K_1}(e_Ne_K) 
    = E^{M_1}_{K_1}(e_N), 
  $
  by \cite{PP:entropy},
  we have 
  \begin{align*}
    E^{M_1}_{K_1}(e_N) 
    &= [M:K]E^{K_1}_M(E^{M_1}_{K_1}(e_N)e_K )e_K \\
    &= [M:K]E^{M_1}_M(e_N)e_K \\ 
    &= \frac{[M:K]}{[M:N]}e_K  \\
    &= \frac{1}{[K:N]}e_K.
  \end{align*}
  Therefore we have
  \begin{displaymath}
    e_Me_Ne_{K_1} = \frac{1}{[K:N]}e_Me_K.
  \end{displaymath}
  And hence we have 
  \begin{displaymath}
    E^{M_2}_{M_1}(e_Me_Ne_{K_1}) = \frac{1}{[K:N]}E^{M_2}_{M_1}(e_M)e_K 
    = \frac{1}{[K:N][M:N]}e_K.  
  \end{displaymath}
  We get the result.
\end{pf}
\begin{lem}\label{lem:isom}
  Let $K$ be an intermediate subfactor for $N \subset M$.
  Let 
  $K \subset M \subset K_1 \subset K_2$ and 
  $N \subset M \subset M_1 \subset M_2$ be the Jones towers 
  for $K \subset M$ and $N \subset M$, respectively.
  If $e_{K_1}$ is the Jones projection for $K_1 \subset M_1$, 
  then there exists a $*$-isomorphism $\varphi$ of $K_2$ onto 
  $e_{K_1}M_2e_{K_1}$ such that $\varphi(x) = xe_{K_1}$ 
  for $x \in K_1$ and 
  $\varphi (e^{K_1}_M) = e_M$,
  where $e^{K_1}_M$ and $e_M$ are the Jones projections for 
  $M \subset K_1$ and $M \subset M_1$, respectively.
\end{lem}
\begin{pf}
  Since $e_{K_1} \in K_1' \subset M'$, it is obvious that 
  $(M \subset K_1) \simeq (Me_{K_1} \subset K_1e_{K_1})$.
  Therefore it is enough to show that $e_{K_1}M_2e_{K_1}$ 
  is the 
  basic extension for $Me_{K_1} \subset K_1e_{K_1}$ with the Jones 
  projection $e_M$.
  By the fact that $e_M = e_{K_1}e_Me_{K_1}$, 
  $e_M$ is an element of $e_{K_1}M_2e_{K_1}$.
  Let $\widetilde{K_2}$ be the basic extension for $K_1 \subset M_1$.
  Since $e_{K_1}\widetilde{K_2}e_{K_1} = K_1e_{K_1}$, 
  we get by Lemma \ref{lem:projection},
  \begin{align*}
    E^{e_{K_1} M_2 e_{K_1}}_{K_1e_{K_1}}(e_M) 
    &= E^{e_{K_1} M_2 e_{K_1}}_{e_{K_1}
      \widetilde{K_2}e_{K_1}}(e_M) \\
    &= e_{K_1} E^{M_2}_{\widetilde{K_2}}(e_M)e_{K_1} \\
    &= E^{M_2}_{\widetilde{K_2}}(e_M) \\
    &= \frac{1}{[M:K]}e_{K_1}.
  \end{align*}
  We can see that 
  \begin{displaymath}
    Me_{K_1} = (K_1 \cap \{e_M\}')e_{K_1} 
    = K_1e_{K_1} \cap \{e_M\}'.
  \end{displaymath}
  Therefore $e_{K_1}M_2e_{K_1}$ is the basic extension for 
  $Me_{K_1} \subset K_1e_{K_1}$ by \cite{PP:entropy}.
\end{pf}

\subsection{Finite dimensional Hopf $C^*$-algebras}

In this subsection we recall 
some facts about finite 
dimensional Hopf $C^*$-algebras.

Let $\H$ be a finite dimensional Hopf $C^*$-algebra 
with a comultiplication $\bigtriangleup_{\H}$
and an anti-pode $S_{\H}$.
Let $\K$ be a subHopf $C^*$-algebra of $\H$, i.e.,
$\K$ is a $*$-subalgebra of $\H$, 
$S_{\H}(\K) \subset \K$ and 
$\bigtriangleup_{\H}(\K) \subset \K \otimes \K$.

\begin{lem}\label{lem:ideal}
Define the subset $\K^{\bot}$ of $\H^*$ by 
$$
\K^{\bot} = \{ \ f \in \H^* \ | (f, k) = 0, \ \forall k \in \K \ \},
$$
where $(\ , \ ): \H^* \times \H \to \Bbb C$ is the dual pairing
defined by $(f, h) = f(h), f \in {\H}^*, h \in \H$.
Then $\K^{\bot}$ is an ideal of $\H^*$.
\end{lem}
\begin{pf}
Let $g$ be an element of $\K^{\bot}$ and $f$ an element of $\H^*$.
Then the element $gf$ of $\H^*$ is determined by the equation 
$$
(gf, h) = (g \otimes f, \bigtriangleup_{\H}(h)), \  \forall h \in \H.
$$
By virtue of  $\bigtriangleup_{\H}(\K) \subset \K \otimes \K$, 
we get
$$
(gf, k) = (g \otimes f, \bigtriangleup_{\H}(k)) = 0, \  \forall k \in \K.
$$
Therefore $gf$ is an element of $\K^{\bot}$.
Similarly, $fg \in \K^{\bot}$.
\end{pf}
By the above lemma, there exists the central projection $p \in \H^*$
such that $\K^{\bot} = p\H^*$.
We put $e_{\K} = 1 - p$.
\begin{prop}
With the above notation, 
the reduced algebra $e_{\K} \H^*$ is the dual 
Hopf $C^*$-algebra of $\K$.
\end{prop}

\begin{pf}
Suppose that $k \in \K$ and $(y,k) = 0, \ \forall y \in e_{\K}\H^*$.
Then 
$$
(f , k) = (e_Kf, k) + (pf, k) = (e_Kf, k) = 0, \ \forall f \in \H^*.
$$
Therefore $k = 0$.
Conversely, suppose that $y \in e_{\K} {\H}^*$ and 
$(y,k) = 0, \ \forall k \in \K$.
Then $y \in K^{\bot} \cap e_{\K}{\H}^* = \{0\}$.
Hence
the form 
$( \ , \ ) |_{e_{\K}{\H}^* \times \K}$ 
establishes a duality between 
$\K$ and $e_{\K}\H^*$. 
So we can identify $e_{\K}{\H}^*$ with  $\K^*$.
Then for $y \in \K^*$ and $k_1, k_2 \in \K$, 
we have
\begin{align*}
(y, k_1k_2) 
&= (\bigtriangleup_{\H^*}(y), k_1 \otimes k_2) \\
&= (\bigtriangleup_{\H^*}(y)(e_{\K} \otimes e_{\K}), k_1\otimes k_2).
\end{align*}
Hence $\bigtriangleup_{{\K}^*}(y) 
= \bigtriangleup_{\H^*}(y)(e_{\K} \otimes e_{\K})$. 
Similarly, we have $S_{\K^*} = S_{\H^*}|_{\K^*}$ 
by the fact that
$$
\overline{(y^*, k^*)} = (S_{\H^*}(y), k), \ 
\forall y \in \K^*, \ \forall k \in \K.
$$
Therefore $e_{\K}{\H}^*$ is again a Hopf $C^*$-algebra with the dual 
algebra $\K$.
\end{pf}

\begin{thm}
Let $\H$ be a finite dimensional Hopf $C^*$-algebra.
The number of subHopf $C^*$-algebras of $\H$ is finite.
\end{thm}
\begin{pf}
By the above proposition, the map $\K \mapsto e_{\K}$ from
the set of  
subHopf $C^*$-algebras of $\H$ 
to central projections of $\H^*$ is injection.
Since the number of central projections of $\H^*$ is finite,
so is that of subHopf $C^*$ algebras of $\H$.
\end{pf}
\begin{Rem}
  Since every finite dimensional Hopf $C^*$-algebra 
  (Kac algebra) admits an \lq\lq outer" action on the 
  AFD II$_1$ factor \cite{Yuchi:Kac},  
  the above theorem immediately follows from
  \cite[Theorem 2.2]{Wa:lattice}.
\end{Rem}

\begin{Defn}
  \label{defn:defHopf}
  Let $\H$ be any Hopf algebra.
  \begin{enumerate}
  \item The left adjoint action of $\H$ on itself is given by 
    \begin{displaymath}
      (ad_lh)(k) = \sum h_1k (S_{\H}(h_2)),
    \end{displaymath}
    for all $h, k \in \H$.
  \item The right adjoint action of $\H$ on itself is given by 
    \begin{displaymath}
      (ad_rh)(k) = \sum (S_{\H}(h_1))kh_2,
    \end{displaymath}
    for all $h, k \in \H$.
  \item A subHopf algebra $\K$ of $\H$ is called 
    {\it normal\ } if both
      \begin{displaymath}
        (ad_l\H)(\K) \subset \K \ \text{and} \ 
        (ad_r\H)(\K) \subset \K.
      \end{displaymath}
  \end{enumerate}
\end{Defn}
\noindent
See \cite[pp.\ 33]{mont:Hopf}.

The next proposition is useful later.
\begin{prop}\label{prop:normalHopf}
  Let $\H$ be  a finite dimensional Hopf algebra with 
  a counit $\varepsilon_{\H}$ and 
  $\K$ a subHopf algebra of $\H$.
  Then $\K$ is normal if and only if 
  $\H\K^+ = \K^+\H$,
  where $\K^+ = \K \cap \ker \varepsilon_{\H}$.
\end{prop}
\noindent
See for a proof \cite[pp.\ 35]{mont:Hopf}.

\subsection{Bimodules }
\label{sec:bimodule}

In this subsection we recall some facts about the 
bimodule calculus associated with  an inclusion of type 
II$_1$ factors (see for example \cite{Oc:QS},\cite{Ya:note}).

Let $A, B, C$ be type II$_1$ factors and let
$
\alpha = \sideset{_A}{_B}H, \beta = 
\sideset{_A}{_B}K, \gamma = \sideset{_B}{_C}L  
$ 
be $A$-$B$, $A$-$B$ and $B$-$C$ Hilbert bimodules, 
respectively.
We write $\alpha \gamma $ for the $A$-$C$ Hilbert bimodule 
$\sideset{_A}{_B}H\otimes_B \sideset{_B}{_C}L$.
We denote by $\< \alpha, \beta \>$ the dimension of the 
space of $A$-$B$ intertwiners from 
$
    \sideset{_A}{_B}H 
$
 to 
$
   \sideset{_A}{_B}K.  
$
The conjugate Hilbert space $H^*$ of $\sideset{_A}{_B}H$ is 
naturally a $B$-$A$ bimodule with $B$-$A$ actions defined by 
\begin{displaymath}
  b \cdot \xi^*\cdot a = (a^* \xi b^*)^*
  \quad \text{for } a \in A \ \text{ and } \  b \in B,
\end{displaymath}
where $\xi^* = \<\cdot, \xi \>_H \in H^*$ for $\xi 
\in \sideset{_A}{_B}H$.
We denote by $\overline{\alpha}$ the conjugate $B$-$A$ 
Hilbert bimodule associated with $\alpha$.
\begin{prop}[Frobenius reciprocity]
  Let $A, B, C$  be type II$_1$ factors, and 
  $\alpha = \sideset{_A}{_B}H, 
  \beta = \sideset{_B}{_C}K $ and 
  $\gamma = \sideset{_A}{_C}L$ be 
  Hilbert bimodules . Then 
  \begin{displaymath}
    \<\alpha \beta, \gamma \> = 
    \<\alpha, \gamma \overline{\beta} \> = 
    \<\beta, \overline{\alpha} \gamma \>.
  \end{displaymath}
\end{prop}
\noindent
See for a proof \cite{Oc:QS}, \cite{Ya:note}.

\begin{exmp}\label{exmp:aut}
  Let $M$ be a type II$_1$ factor with the normalized 
  trace $\tau_M$. 
  As usual we let $L^2(M)$ be the Hilbert space obtained by 
  completing $M$ in the norm $\parallel x \parallel_2 
  = \sqrt{\tau_M(x^*x)}, \ x \in M$. 
  Let $\eta: M \to L^2(M)$ be the canonical implementation.
  Let $J : L^2(M) \to L^2(M)$ be the modular conjugation 
  defined by $J\eta(x) = x^*, \  x \in M$.
  For $\theta \in \Aut(M)$, we define 
  $\sideset{_M}{_M}{L^2(\theta)}$,
  the $M$-$M$ Hilbert bimodule, 
  by
  \begin{enumerate}
  \item $\sideset{_M}{_M}{L^2(\theta)} = L^2(M)$ as 
    a Hilbert space,
  \item $x \cdot \xi \cdot y = x J \theta(y)^*J \xi, \  
    x, y \in M, \ \xi \in L^2(M)$.
  \end{enumerate}
  Then for $\theta, \theta_1, \theta_2 \in \Aut(M)$ 
  we have 
  \begin{displaymath}
    \overline{\sideset{_M}{_M}{L^2(\theta)}}
    \simeq \sideset{_M}{_M}{L^2(\theta^{-1})} 
  \end{displaymath}
  \begin{displaymath}
    \sideset{_M}{_M}{L^2(\theta_1)} 
    \underset{M}\otimes 
    \sideset{_M}{_M}{L^2(\theta_2)} 
    \simeq  \sideset{_M}{_M}{L^2(\theta_1\theta_2)}. 
  \end{displaymath}
\end{exmp}

A bimodule $\alpha = \sideset{_A}{_B}H$ is called 
irreducible if $\<\alpha, \alpha\> = 1$,
 i.e., $\oEnd_{\AB}(\sideset{_A}{_B}H) \simeq \Bbb C$. 
If $\<\alpha, \alpha\> < \infty$,
$\alpha = \sideset{_A}{_B}H$, 
then we can get 
an $A$-$B$ irreducible bimodule by cutting $\sideset{_A}{_B}H$ 
by a minimal projection in $\oEnd_{\AB}(\sideset{_A}{_B}H)$.

\begin{exmp}\label{exmp:bimodule}
  Let $N \subset M$ be an inclusion of type II$_1$ factors.
  We define the $N$-$M$ bimodule $\sideset{_N}{_M}{L^2(M)}$by
  actions
  \begin{displaymath}
    x \cdot \xi \cdot y = xJy^*J\xi,\quad \xi \in L^2(M),  x \in N, y \in M.
  \end{displaymath}
  Then we can see that 
  \begin{displaymath}
    \operatorname{End}(\sideset{_N}{_M}{L^2(M)}) \simeq N' \cap M. 
  \end{displaymath}
  In particular, 
  if $N' \cap M = \Bbb C$, then $\sideset{_N}{_M}{L^2(M)}$ is an 
  irreducible $N$-$M$ bimodule.
\end{exmp}
The next lemma is well known.
\begin{lem}\label{lem:equivalent}
  Let $N \subset M$ be a pair of type II$_1$ factors with 
  finite index and $M_1$ the basic extension for the inclusion 
  $N \subset M$.
  For $\theta \in \Aut(N)$, 
  $\sideset{_N}{_{\theta(N)}}{L^2(M)} \simeq   
  \sideset{_N}{_N}{L^2(M)}$ if and only if there exists a 
  unitary $u \in M_1$ such that $uxu^* = \theta(x)$, for all 
  $x \in N$, where $\sideset{_N}{_{\theta(N)}}{L^2(M)}$
  is defined as in Example \ref{exmp:bimodule}.
\end{lem}
%\begin{pf}
%  Let $T$ be an onto isometry from 
%  $\sideset{_N}{_{\theta(N)}}{L^2(M)}$ to 
%  $\sideset{_N}{_N}{L^2(M)}$ such that 
%  \begin{displaymath}
%    T(x \cdot \xi \cdot y) = x \cdot T(\xi) \cdot y
%  \end{displaymath}
%  for $\forall \xi$ and $x, y \in N$.
%  Since $T(x\xi) = xT(\xi)$ for  
%  $\forall \xi \in L^2(M)$ and $\forall x \in N$,
%  $T \in N'$ on $L^2(M)$.
%  Since $T(J\theta(y)^*J\xi) = Jy^*JT (\xi)$, 
%  $\forall \xi \in L^2(M)$ and $\forall y \in N$,
%  we have 
%  \begin{displaymath}
%    TJ\theta(x)J = JxJT, \ \forall x \in N,
%  \end{displaymath}
%  where $J$ is the modular conjugation on $L^2(M)$.
%  Put $ u = JT^*J \in JN'J   (= M_1)$.
%  Then $uxu^* = \theta(x)$, $\forall x \in N$.
%  Conversely, suppose that there exists a 
%  unitary $u \in M_1$ such that 
%  $uxu^* = \theta(x)$, $ \forall x \in N$.
%  Define $\phi: \sideset{_N}{_{\theta(N)}}{L^2(M)}
%  \to \sideset{_N}{_N}{L^2(M)}$ by 
%  $\phi(\xi) = Ju^*J\xi$, $\xi \in L^2(M)$.
%  Then since $Ju^*J \in N'$,
%  $\phi(x\xi) = x \phi(\xi) $, $\forall x \in N$.
%  By the definition of $\phi$, we have 
%  $\phi(J\theta(y)^*J\xi) = Jy^*u^*J\xi 
%  = Jy^*J \phi(\xi)$, $\forall y \in N$.
%  We have thus proved the lemma.    
%\end{pf}

\begin{exmp}\label{exmp:group}
  Let $\gamma : G \to \Aut(N)$ be an outer action of  
  a finite group $G$ on a type II$_1$ factor $N$.
  Let $M = N \rtimes_{\gamma}G$ 
  be the crossed product 
  and $\rho$ the $N$-$M$ bimodule $\sideset{_N}{_M}{L^2(M)}$ 
  defined as in Example \ref{exmp:bimodule}. 
  If $\{ \lambda_g | \ g \in G \} $ is a unitary implementation
  for the crossed product, then 
  each element $x \in M$ is written in the form 
  $x = \sum_{g \in G}x_g\lambda_g, x_g \in N$.
  This implies that the irreducible decomposition of 
  $\rho\overline{\rho} = \sideset{_N}{_N}{L^2(M)}$ is
  \begin{displaymath}
    \bigoplus_{g \in G} 
    \sideset{_N}{_N}{(\overline{N\lambda_g}^{\Vert \cdot \Vert_2})} 
    \simeq
    \bigoplus_{g \in G} 
    \sideset{_N}{_N}{L^2(\gamma_g)},
  \end{displaymath}
  where $\sideset{_N}{_N}{L^2(\gamma_g)}$ is the 
  $N$-$N$ bimodule as in Example \ref{exmp:aut}.
\end{exmp}

\section{Definition of normal intermediate subfactors}
\label{sec:def}

In this section, we shall introduce the notion of normality for 
intermediate subfactors and study its properties.

Let $N \subset M$ be a pair of type II$_1$ factors 
with $[M:N] < \infty$. 
Let 
$
N \subset M \subset M_1 \subset M_2
$
be the Jones tower of $N \subset M$, obtained by 
iterating the basic extensions.
We denote by $\cal F$, Ocneanu's Fourier transform from 
$N' \cap M_1$ onto $M' \cap M_2$ 
i.e., 
\begin{displaymath}
  \cal F (x) 
  = [M:N]^{- \frac{3}{2}}E^{N'}_{M'}(xe_Ne_M), 
  \  x \in N' \cap M_1, 
\end{displaymath}
where $E^{N'}_{M'}$ is the conditional expectation 
from $N'$ onto $M'$.
\begin{defn}
  Let $K$ be an intermediate subfactor of $N \subset M$ 
  and $e_K$ the Jones projection for the inclusion 
  $K \subset M$. 
  Then we call that $K$ is {\it normal\ } in $N \subset M$
  if $e_K$ and $\cal F (e_K)$ are elements of the centers 
  of $N' \cap M_1$ and $M' \cap M_2$, respectively.
\end{defn}
\begin{lem}\label{lem:dual}
  Let $K$ be an intermediate subfactor for an irreducible 
  inclusion $N \subset M$ of type II$_1$ factors 
  with finite index.
  Let $K_1$ and $M_1$ be the basic extensions 
  for $K \subset M$ and $N \subset M$, respectively.
  Then $K$ is normal in $N \subset M$ if and only if 
  $K_1$ is normal in $M \subset M_1$.
\end{lem}
\begin{pf}
  Since $\cal F(e_K) = \lambda e_{K_1}$
  for some $\lambda \in \Bbb C$, 
  It is obvious by the definition.
\end{pf}
\begin{prop}\label{prop:fix}
  Let $N$ be the fixed point algebra $M^{(G, \alpha)}$   
  of a type II$_1$ factor $M$ by an outer action 
  $\alpha$ of a finite 
  group $G$. If $K = M^{(H, \alpha)}$ is an  
  intermediate subfactor  
  associated with 
  a subgroup $H$ of $G$, 
  then $K$ is normal in $N \subset M$ if and only if 
  $H$ is a normal subgroup of $G$.
\end{prop}
\begin{pf}
  Let $\{u_g | \ g \in G\}$ be unitary operators on $L^2(M)$ 
  defined by $u_g\eta(x) = \eta(\alpha_g(x))$, $x \in M$, 
  where $L^2(M)$ and $\eta$ are defined as in Example \ref{exmp:aut}.
  Then $N = M \cap \{u_g | \ g \in G\}'$, 
  $M_1 = (M \cup \{u_g | \ g \in G\})''$ 
  and $N' \cap M_1 = \{u_g | \ g \in G\}'' \simeq \Bbb C G$.
  Since $K = M^H = M \cap \{ u_h | \ \in H \}'$, 
  the Jones projection $e_K$ for $K \subset M$ 
  is 
  $
    \frac{1}{\sideset{^\#}{}H}\sum_{h \in H}
    u_h.
  $
  Since 
  \begin{displaymath}
    u_g e_K u_g^* = \frac{1}{\sideset{^\#}{}H}\sum_{h \in H} u_{ghg^{-1}}
  \quad \text{for any } \ g \in G,
  \end{displaymath}
  $H$ is a normal subgroup if and only if 
  $e_K$ is an element of the center of $N'\cap M_1$.
  Since $M' \cap M_2$ is a commutative algebra, 
  $\cal F(e_K)$ is always an element of the center of 
  $M'\cap M_2$. So we get the result.  
\end{pf}
\begin{prop}\label{prop:crossed}
  Let $M$ be the crossed product $N \rtimes_{\alpha}G$ 
  of a II$_1$ factor $N$ by an outer action $\alpha$ of a finite 
  group $G$. If $K = N \rtimes_{\alpha}H$ 
  is an intermediate subfactor associated with a 
  subgroup $H$ of $G$, then $K$ is normal in $N \subset M$ i
  f and only if 
  $H$ is a normal subgroup of $G$.
\end{prop}
\begin{pf}
  This immediately follows from  Lemma \ref{lem:dual}  
  and Proposition \ref{prop:fix}.
\end{pf}
\begin{prop}
  Let $\alpha : G \to \Aut(P)$ be an outer action of 
  a finite group $G$ on a type II$_1$ factor $P$ and 
  $H$ a subgroup of $G$.
  Let $M$ be the fixed point algebra 
  $P^{(H, \alpha)}$ and 
  $N$ the fixed point algebra  $P^{(G, \alpha)}$.
  For $K \in \Lat(N \subset M)$, there is a  
  subgroup $A$ of $G$ such that $H \subset A \subset G$ and 
  $K = P^{(A, \alpha)}$. Then 
  $K$ is a normal intermediate subfactor of $N \subset M$ 
  if and only if $AgH = HgA$, for $\forall g \in G$.
\end{prop}
\begin{pf}
  Let $\{u_g | \ g \in G\}$ be unitary operators on $L^2(P)$ 
  defined by $u_g\eta(x) = \eta(\alpha_g(x))$, $x \in P$, 
  where $L^2(P)$ and $\eta$ are defined as in Example \ref{exmp:aut}.
  Let $P_1$ be the basic extension for $N \subset P$.
  Then 
  \begin{displaymath}
    N' \cap P_1 
    = \{\sum_{g \in G}x_gu_g \ | \ x_g \in \Bbb C \} 
    \simeq \Bbb C G.
  \end{displaymath}
  Let $e^P_M$ be the Jones projection for 
  $M \subset P$. Then 
  \begin{displaymath}
    e^P_M = \frac{1}{\sideset{^\#}{}H}\sum_{h \in H}
    u_h.
  \end{displaymath}
  Let $M_1$ be the basic extension for $N \subset M$.
  Then by Lemma \ref{lem:isom}, 
  \begin{align*}
    N' \cap M_1 
    &\simeq e^P_M(N' \cap P_1)e^P_M \\
    &= \{ \sum_{g \in G}  \sum_{h, k \in H} x_gu_{hgk} \ 
    | \ x_g \in \Bbb C \}.
  \end{align*}
  Therefore
  \begin{align*}
    e^M_K \in \cal Z(N' \cap M_1) 
    &\Leftrightarrow e^P_Me^P_Ke^P_M (= e^P_K =  
    \frac{1}{\sideset{^\#}{}A}\sum_{a \in A}u_a)
    \in \cal Z (e^P_M(N' \cap P_1)e^P_M) \\
    &\Leftrightarrow \sum_{a \in A}\sum_{h, k \in H}
    u_{ahgk}  = \sum_{a \in A}\sum_{h, k \in H}u_{hgka} 
    \ \text{for } \forall g \in G \\
    &\Leftrightarrow AgH = HgA \ \text{for } \forall g \in G.
  \end{align*}
  Since $M'\cap M_2$ is a commutative algebra, we get the 
  the result.
\end{pf}
\begin{prop}
  Let $\alpha : G \to \Aut(P)$ be an outer action of 
  a finite group $G$ on a type II$_1$ factor $P$ and 
  $H$ a subgroup of $G$.
  Let $M$ be the crossed product 
  $P \rtimes_{\alpha}G$ and 
  $N$ the crossed product $P \rtimes_{\alpha}H$.
  For $K \in \Lat(N \subset M)$, there is a  
  subgroup $A$ of $G$ such that $H \subset A \subset G$ and 
  $K = P \rtimes_{\alpha} A$. Then 
  $K$ is a normal intermediate subfactor of $N \subset M$ 
  if and only if $AgH = HgA$, for $\forall g \in G$.
\end{prop}
\begin{pf}
  This immediately follows from Lemma \ref{lem:dual} 
  and the above proposition.
\end{pf}
\begin{prop}
  Let $N \subset M$ and $Q \subset P$ 
  be irreducible inclusions 
  of type II$_1$ factors
  with finite indices.
  Then both of $N \otimes P$ and $M \otimes Q$ 
  are normal intermediate subfactors of
  $N \otimes Q \subset M \otimes P$.
\end{prop}
\begin{pf}
  Let $M_1 = \<M,e_N\>$ and 
  $P_1 = \<P, e_Q\>$ be the basic extension for 
  $N \subset M$ and $Q \subset P$ with the 
  Jones projections $e_N$ and $e_Q$, respectively.
  Then $M_1 \otimes P_1$ is the basic extension for 
  $N \otimes Q \subset M \otimes P$ with the 
  Jones projection $e_N \otimes e_Q$.
  Moreover, $e_N \otimes 1$ and $1 \otimes e_Q$ are 
  the Jones projections for $N \otimes P \subset M \otimes P$ 
  and $M \otimes Q \subset M \otimes P$, respectively.
  Since  $N \subset M$ and $Q \subset P$ are irreducible, 
  $e_N$ and $e_Q$ are elements of the centers of 
  $N' \cap M_1$ and $Q' \cap P_1$, respectively 
  by \cite[Proposition 1.9]{PP:entropy}.
  And hence $e_N \otimes 1$ and $1 \otimes e_Q$ are 
  elements of the center of $(N \otimes Q)' \cap (M_1 \otimes P_1) 
  = (N'\cap M_1) \otimes (Q' \cap P_1)$.
  Similarly, we can observe that 
  $\cal F(e_N \otimes 1)$ and $\cal F(1 \otimes e_Q)$ are 
  elements of the center of 
  $(M \otimes P)' \cap (M_2 \otimes  P_2)$, 
  where $M_2$ and $P_2$ are the basic extension for 
  $M \subset M_1$ and $P \subset P_1$, respectively.
  We have thus proved the proposition
\end{pf}

In \cite{Wa:lattice} Y.~Watatani introduced the notion of 
{\it quasi-normal
intermediate subfactors\/ } 
to study the modular identity for intermediate subfactor
lattices.

\begin{Defn}
  Let $N \subset M$ be an inclusion of type II$_1$ 
  factors with finite index and $K$ an intermediate 
  subfactor of $N \subset M$.
  Then $K$ is {\it quasi-normal\/} 
  (or {\it doubly commuting \/}) 
  if for any $L \in \Lat(N \subset M)$,
  \begin{alignat*}{3}
         K &  \ &  \subset  & \ & \quad K  & \lor  L  \\
       \cup&  \ &     \   &  \   &  \ & \cup  \  \\
    K \land&  L \quad  &  \subset & \ &  \ & L  \   
  \end{alignat*}
  and
  \begin{alignat*}{3}
         K_1 &  \ &  \subset  & \ & \quad K_1  & \lor  L_1  \\
       \cup&  \ &     \   &  \   &  \ & \cup  \  \\
    K_1 \land&  L_1 \quad  &  \subset & \ &  \ & L_1  \   
  \end{alignat*}
  are commuting squares  (see for example \cite{GHJ:coxeter}), 
  where $K_1$ and $L_1$ are the basic 
  extension for $K \subset M$ and $L \subset M$, respectively.
\end{Defn}
\begin{prop}\label{prop:quasi}
  Let $N \subset M$ be an irreducible inclusion of type II$_1$ 
  factors with finite index.
  If $K$ is a normal intermediate subfactor of  $N \subset
  M$ then $K$ is quasi-normal in $N \subset M$
\end{prop}
\begin{pf}
  Suppose that the Jones projection $e_K$ for 
  $K \subset M$ is an element of the center of $N' \cap M_1$. 
  Then since
  for any intermediate subfactor $L$ of $N \subset M$, 
  the Jones projection $e^{K \vee L}_K$ for 
  $K \subset (K \vee L)$ is also a central
  projection in $K' \cap (K \vee L)_1$,   
  we have 
  \begin{alignat*}{3}
    K &  \ &  \subset  & \ & \quad K  & \lor  L  \\
    \cup&  \ &     \   &  \   &  \ & \cup  \  \\
    K \land&  L \quad  &  \subset & \ &  \ & L  \   
  \end{alignat*}
  is a commuting square.
  Similarly, if $\cal F (e_K)$ is an element of the center of 
  $M' \cap M_2$, then 
  \begin{alignat*}{3}
    K_1 &  \ &  \subset  & \ & \quad K_1  & \lor  L_1  \\
    \cup&  \ &     \   &  \   &  \ & \cup  \  \\
    K_1 \land&  L_1 \quad  &  \subset & \ &  \ & L_1  \   
  \end{alignat*}
  is a commuting square.
  Therefore if $K$ is normal in $N \subset M$, then 
  $K$ is quasi-normal.
\end{pf}

We have a characterization of normal intermediate 
subfactors in terms of bimodules.
Let $K$ be an intermediate subfactor of an irreducible
inclusion $N \subset M$ of type II$_1$ factors with finite 
index. We note that $e_K$ is in the center of $N' \cap M_1$ 
if and only if for any $T \in \oEnd(\sideset{_N}{_N}{L^2(M)})$, 
$TL^2(K) \subset L^2(K)$.

\begin{prop}\label{lem:defbi}
  Let $K$ be an intermediate subfactor for an irreducible 
  inclusion 
  $N \subset M$ of type II$_1$ factors with finite index.
  Let $\alpha$ be the $N$-$K$ bimodule $\sideset{_N}{_K}{L^2(K)}$ 
  and $\beta$ the $K$-$M$ bimodule $\sideset{_K}{_M}{L^2(M)}$.
  If $\rho$ is the $N$-$M$ bimodule 
  $\alpha \beta = \sideset{_N}{_M}{L^2(M)}$,
  then  $K$ is normal in $N \subset M$ if and only if 
  \begin{enumerate}
  \item $\<\alpha \overline{\alpha}, \rho \overline{\rho} \> 
    = \<\alpha \overline{\alpha}, \alpha \overline{\alpha}\>$, 
  \item $\<\overline{\beta}\beta, \overline{\rho}\rho\> 
    = \<\overline{\beta}\beta, \overline{\beta}\beta \>$.
  \end{enumerate}
\end{prop}
\begin{pf}
  Since $ \oEnd( \sideset{_N}{_K}{L^2(K)} ) =
  N' \cap \<N, e^K_N\> \simeq 
  e_K(N' \cap M_1)e_K $ by Lemma \ref{lem:isom}, 
  if $e_K$ is an element of the center of $N' \cap M_1$, 
  then for any irreducible $N$-$N$ bimodule $\sigma$ 
  contained in $\alpha\overline{\alpha}$, 
  the multiplicity of $\sigma$ in $\alpha\overline{\alpha}$ 
  is equal to the multiplicity of $\sigma$ in 
  $\rho \overline{\rho}$.
  Therefore we have 
  \begin{displaymath}
    \<\alpha \overline{\alpha}, \rho \overline{\rho} \> 
    = \<\alpha \overline{\alpha}, \alpha \overline{\alpha}\>.
  \end{displaymath}
  Conversely, suppose that 
  $e_K$ is not an element of the center of $N' \cap M_1$. 
  Then there exist minimal projections $p, q \in 
  e_K(N' \cap M_1)e_K$ 
  such that
  \begin{displaymath}
    p \sim q \ \text{in} \ (N' \cap M_1) \quad 
    \text{and}\quad p \not\sim q \ \text{in} \  
    e_K(N' \cap M_1)e_K.
  \end{displaymath}
  Therefore we have 
  \begin{displaymath}
    \<\alpha \overline{ \alpha}, \rho \overline{\rho} \> 
    \not= \<\alpha \overline{\alpha}, \alpha \overline{\alpha}\>.
  \end{displaymath}
  And hence $e_K$ is an element of the center of $(N' \cap M_1) $
  if and only if 
  \begin{displaymath}
    \<\alpha \overline{\alpha}, \rho \overline{\rho} \> 
    = \<\alpha \overline{\alpha}, \alpha \overline{\alpha}\>.
  \end{displaymath}
  Similarly, we can see that 
  $e_{K_1}$ is an element of the center of $(M' \cap M_2)$ 
  if and only if 
  \begin{displaymath}
    \<\overline{\beta}\beta, \overline{\rho}\rho \> 
    =  \<\overline{\beta}\beta,\overline{\beta}\beta \>.
  \end{displaymath}
  Since $\cal F(e_K) = \lambda e_{K_1}$ for some  
  $\lambda \in \Bbb C$, 
  we get the result.
\end{pf}

\begin{thm}\label{thm:depth2}
  Let $K$ be an intermediate subfactor for an irreducible
  inclusion $N \subset M$ of type II$_1$ factors 
  with finite index. If the depths of
  $N \subset K$ and $K \subset M$ are both 2, then $K$ is
  normal in $N \subset M$. 
\end{thm}

\begin{pf}
  Let $\alpha$ be the $N$-$K$ bimodule $\sideset{_N}{_K}{L^2(K)}$ 
 and $\beta$ the $K$-$M$ bimodule $\sideset{_K}{_M}{L^2(M)}$.
 By the assumption, we have 
 \begin{displaymath}
   \alpha \overline{\alpha}\alpha \simeq \underbrace{\alpha \oplus 
   \alpha \oplus \cdots \oplus \alpha}_{[K:N] \text{times}}  
 \end{displaymath}
and
\begin{displaymath}
  \overline{\beta}\beta \overline{\beta} \simeq
  \underbrace{\overline{\beta} \oplus \overline{\beta} \oplus
    \cdots \oplus \overline{\beta}}_{[M:K] \text{times}}.
\end{displaymath}
And hence
\begin{displaymath}
  \<\alpha\overline{\alpha}, \alpha\overline{\alpha} \> = \<
  \alpha \overline{\alpha} \alpha, \alpha \> = [K:N]
\end{displaymath}
and 
\begin{displaymath}
  \<\overline{\beta}\beta, \overline{\beta}\beta \> 
  = \< \overline{\beta}\beta\overline{\beta}, \overline{\beta}\> 
  = [M:K]
\end{displaymath}
by Frobenius reciprocity.
Since $N \subset M$ is irreducible,  if $\rho$ is the $N$-$M$ bimodule
$\sideset{_N}{_M}{L^2(M)} (= \alpha\beta)$, then 
\begin{displaymath}
  1 = \< \rho, \rho\> = \<\alpha\beta, \alpha\beta \> = \<
  \overline{\alpha} \alpha, \beta \overline{\beta} \>.
\end{displaymath}
And hence we have  
\begin{align*}
  \< \alpha \overline{\alpha}, \rho \overline{\rho} \> 
  &= \< \alpha \overline{\alpha}, 
  \alpha \beta \overline{\beta} \overline{\alpha} \> \\
  &= \< \alpha \overline{\alpha} \alpha, 
  \alpha \beta \overline{\beta} \> \\
  &= [K:N] \< \alpha, \alpha \beta \overline{\beta} \> \\
  &= [K:N] \< \overline{\alpha}\alpha, 
  \beta  \overline{\beta} \>  = [K:N], 
\end{align*}
i.e., 
\begin{displaymath}
  \<\alpha\overline{\alpha}, \rho\overline{\rho} \> =
  \<\alpha\overline{\alpha}, \alpha\overline{\alpha}\>.
\end{displaymath}
Similarly, we have 
\begin{displaymath}
  \<\overline{\beta}\beta, \overline{\rho}\rho \> = 
  \<\overline{\beta}\beta, \overline{\beta}\beta \>.
\end{displaymath}
So we get the result by Lemma \ref{lem:defbi}.
 \end{pf}

 \begin{prop}\label{prop:int}
   Let $M_0, N_0, K$ be intermediate subfactors for an 
   irreducible inclusion $N \subset M$ of type II$_1$ 
   factors with finite index such that 
   \begin{displaymath}
     N \subset N_0 \subset K \subset M_0 \subset M.
   \end{displaymath}
   If $K$ is normal in $N \subset M$, then $K$ is also normal in 
   $N_0 \subset M_0$.
 \end{prop}

 \begin{pf}
   Let $\alpha = \sideset{_N}{_K}{L^2(K)}$, 
   $\alpha_0 = \sideset{_{N_0}}{_K}{L^2(K)}$,  
   $\beta = \sideset{_K}{_M}{L^2(M)}$ 
   and $\beta_0 = \sideset{_K}{_{M_0}}{L^(M_0)}$.
   Since 
   $$
   \alpha \overline{\alpha} = \sideset{_N}{_K}{L^2(K)} 
   \otimes_K \sideset{_K}{_N}{L^2(K)} = \sideset{_N}{_K}{L^2(K)}
   \otimes \sideset{_K}{_K}{L^2(K)} \otimes_K 
   \sideset{_K}{_N}{L^2(K)},   
   $$
   we have 
   \begin{displaymath}
     \< \alpha \overline{\alpha}, \alpha \overline{\alpha} \>
     = \< \overline{\alpha}\alpha \overline{\alpha}\alpha,
     \sideset{_K}{_K}{L^2(K)} \> 
   \end{displaymath}
   by Frobenius reciprocity.
   Since 
   $
   \<\alpha \overline{\alpha}, 
   \alpha\beta \overline{\beta}\overline{\alpha}\> =
   \< \alpha \overline{\alpha},
   \alpha\overline{\alpha} \>
   $
   by the assumption, we have 
   \begin{displaymath}
     \<\overline{\alpha}\alpha \overline{\alpha}\alpha,
     \beta \overline{\beta} \> = 
     \< \overline{\alpha}\alpha \overline{\alpha}\alpha,
     \sideset{_K}{_K}{L^2(K)} \>, 
   \end{displaymath}
   i.e., the irreducible $K$-$K$ sub-bimodules of 
   $\overline{\alpha}\alpha \overline{\alpha}\alpha$ 
   contained in 
   $\beta\overline{\beta}$ is only $\sideset{_K}{_K}{L^2(K)}$.
   Since $\overline{\alpha_0}\alpha_0$ is contained in 
   $\overline{\alpha}\alpha$ and $\beta_0\overline{\beta_0}$ 
   is contained in $\beta\overline{\beta}$, we have 
   \begin{displaymath}
     \<\overline{\alpha_0}\alpha_0 \overline{\alpha_0}\alpha_0,
     \beta_0 \overline{\beta_0} \> 
     = \<\overline{\alpha_0}\alpha_0 \overline{\alpha_0}\alpha_0, 
     \sideset{_K}{_K}{L^2(K)} \>, 
   \end{displaymath}
   i.e., 
   \begin{displaymath}
     \<\alpha_0\overline{\alpha_0}, 
     \alpha_0\beta_0 \overline{\beta_0} \overline{\alpha_0} \>
     = \<\alpha_0\overline{\alpha_0}, 
     \alpha_0\overline{\alpha_0} \>.
   \end{displaymath}
   By the same argument, we have 
   \begin{displaymath}
     \< \overline{\beta_0}\beta_0, 
     \overline{\beta_0} \overline{\alpha_0} \alpha_0\beta_0 \>
     =  \< \overline{\beta_0} \beta_0, 
     \overline{\beta_0} \beta_0 \>.
   \end{displaymath}
   We have thus proved the proposition.
 \end{pf}

\section{Normal intermediate subfactors for depth 2 inclusions}

It is well-known that the crossed product of 
a finite dimensional Hopf $C^*$ algebra 
(Kac algebra) is  characterized by the depth 2 condition.
In this section we study normal intermediate subfactors 
for depth 2 inclusions.

\subsection{The action of $K' \cap K_1$ on $M$}

Let $N \subset M$ be an irreducible, depth 2 inclusion 
of type II$_1$ factors with 
finite index. 
Let $N \subset M \subset M_1 \subset M_2$ 
be the Jones tower for $N\subset M$.
We put $ A = N' \cap M_1$ and $B = M' \cap M_2$.
Then $A$ and $B$ are dual pair 
of Hopf $C^*$-algebras 
with a pairing 
\begin{displaymath}
  (a, b) = [M:N]^2 \tau(a e_Me_N b), \ 
  \text{for} a \in A \ 
  \text{and} \  b \in B, 
\end{displaymath}
where $e_N$ and $e_M$ are the Jones projections for 
$N \subset M$ and $M \subset M_1$, respectively.
Define a bilinear map $A \times M \to M$ 
(denoted by $a \odot x$) by setting
\begin{displaymath}
  a \odot  x = [M:N]E^{M_1}_M(axe_N),
\end{displaymath}
for  $\ x\in M$ and $a \in A$ 
This map is 
a left action of Hopf $C^*$ algebra $A$ and  
\begin{displaymath}
  N = M^A = \{\ x \in M \ | \ a \odot x = \varepsilon(a)x, \
\forall a \in A \}, 
\end{displaymath}
where $\varepsilon : A \to \Bbb C$ is the counit determined by 
$ae_N = \varepsilon(a)e_N$ (see \cite{Szym:Hopf}).

\begin{prop}\label{prop:subalgebra}
  Let $K$ be an intermediate subfactor of $N \subset M$ and 
  $K_1$ the basic extension for $K \subset M$.
  We put $H = K' \cap K_1$.
  If $a$ is an element of $H$, then 
  \begin{displaymath}
    [M:K]E^{K_1}_M(axe_K) = [M:N]E^{M_1}_M(axe_N), \ \forall x \in M.
  \end{displaymath}
  This implies
  \begin{displaymath}
    K = M^H = \{ \ x \in M \ | \ a \odot x = \varepsilon (a)x, 
    \ \forall a \in H \ \}. 
  \end{displaymath}
\end{prop}
\begin{pf}
  Since 
  \begin{math}
    e_K = \frac{[M:N]}{[M:K]}E^{M_1}_{K_1}(e_N) 
  \end{math} 
  by \cite{SW:angle},  
  we have 
  \begin{align*}
    [M:K]E^{K_1}_M(axe_K) 
    &= [M:K]E^{K_1}_M(ax\frac{[M:N]}{[M:K]}E^{M_1}_{K_1}(e_N)) \\
    &= [M:N]E^{K_1}_M(E^{M_1}_{K_1}(axe_N)) \\
    &= [M:N]E^{M_1}_M(axe_N)
  \end{align*}
  for $\forall a \in H$ and $\forall x \in M$.
\end{pf}

\subsection{Hopf algebra structures on $K' \cap K_1$}

Let $N \subset M$ be an irreducible, depth 2 inclusion 
of type II$_1$ factors with 
finite index and 
$K$ an intermediate subfactor of $N \subset M$.
Then the depth of $K \subset M$ is not 2 in general.
In this subsection we shall prove that 
if the depth of $K \subset M$ is 2, then 
$H = K' \cap K_1$ is a subHopf $C^*$ algebra of 
$A = N' \cap M_1$.
 
By Lemma\ref{lem:isom},
there exists an isomorphism $\varphi$ 
of $K_2$ onto $e_{K_1}M_2 e_{K_1}$ such that 
$\varphi(x) = xe_{K_1}$ for $x \in K_1$ and 
$\varphi(e^{K_1}_M) = e_M$, 
where $K \subset M \subset K_1 \subset K_2$ 
is the Jones tower for the inclusion $K \subset M$
and $e^{K_1}_M$ is the Jones projection for $M \subset K_1$. 
\begin{lem}\label{lem:pairing}
  With the above notation, we have 
  \begin{displaymath}
    [M:K]^2 \tau(h e^{K_1}_M e_K k) 
    =  [M:N]^2\tau(h e_M e_N \varphi(k))
  \end{displaymath}
  for $\forall h \in H = K' \cap K_1$
  and $\forall k \in     M' \cap K_2$.
\end{lem}
\begin{pf}
  By the fact that $e_{K_1}e_N e_{K_1} 
  = E^{M_1}_{K_1}(e_N)e_{K_1}
  = \frac{[M:K]}{[M:N]} e_Ke_{K_1}$, 
  we have 
  $\varphi(e_K) = e_Ke_{K_1} = \frac{[M:N]}{[M:K]} e_{K_1}e_N e_{K_1}$.
  Therefore
  \begin{align*}
    [M:K]^2 \tau(h e^{K_1}_M e_K k)
    &= [M:K]^2 [K:N]\tau(\varphi(h e^{K_1}_M e_K k)) \\
    &= [M:K]^2 [K:N] \frac{[M:N]}{[M:K]} 
    \tau(\varphi(h) e_M  e_{K_1}e_N e_{K_1} \varphi(k)) \\
    &=    [M:N]^2\tau(h e_M e_N \varphi(k)).
  \end{align*}
\end{pf}
\begin{lem}\label{lem:A}
  Let $N \subset M$ be an irreducible, depth 2 inclusion
  of type II$_1$ factors with finite index and 
  $K$ an intermediate subfactor for $N \subset M$.
  Let $N \subset M \subset M_1 \subset M_2 $ 
  and $ K \subset M \subset K_1 \subset K_2$ be 
  the Jones towers for  $ N \subset M $ and $ K \subset M$, 
  respectively.
  If the depth of $K \subset M$ is 2, then
  for any $b \in M' \cap M_2$, 
  there exist 
  elements $\{x_i\}, \{y_i\}$ of $N' \cap M_1$ such that 
  \begin{displaymath}
    b = \sum_i x_i e_M y_i
  \end{displaymath}
  and 
  \begin{displaymath}
    \sum_i E^{M_1}_{K_1}(x_i)e_M E^{M_1}_{K_1}(y_i) 
    \in (K'\cap K_1)e_M(K' \cap K_1),
  \end{displaymath}
  where $E^{M_1}_{K_1}$ is 
  the trace preserving conditional expectation from $M_1$ 
  onto $K_1$.
\end{lem}
\begin{pf}
  Since the depth of $N \subset M$ is 2, 
  \begin{displaymath}
    (N' \cap M_1)e_M(N' \cap M_1) = N' \cap M_2.
  \end{displaymath}
  And hence any element $b \in M' \cap M_2$ is written in 
  the form
  \begin{displaymath}
    b = \sum_i x_i e_M y_i, \quad x_i,y_i \in N'\cap M_1.
  \end{displaymath}
  Since the depth of $K \subset M$ is 2,
  \begin{displaymath}
    (K' \cap K_1) e^{K_1}_M(K' \cap K_1) 
   = K' \cap K_2,
 \end{displaymath}
 where $e^{K_1}_M$ is the Jones projection for 
 $M \subset K_1$.
 By Lemma \ref{lem:isom}, we have
 \begin{displaymath}
   (K' \cap K_1) e_M (K' \cap K_1) = e_{K_1}(K' \cap M_2)e_{K_1}. 
 \end{displaymath}
 Therefore we have 
 \begin{align*}
   e_{K_1}be_{K_1} 
   &= e_{K_1}(\sum_i x_i e_M y_i) e_{K_1}\\
   &= \sum_i E^{M_1}_{K_1}(x_i)e_M E^{M_1}_{K_1}(y_i)
   \in (K' \cap K_1) e_M (K' \cap K_1).
 \end{align*}
 we have thus proved the lemma.
\end{pf}
\begin{prop}\label{prop:subHopf}
  Suppose that the depth of $N \subset M$ is 2.
  Let $K$ be an intermediate subfactor for $N \subset M$ and 
  $K_1$ the basic extension for $K \subset M$.
  If the depth of $K \subset M$ is 2, 
  then  $H = K' \cap K_1$ is 
  a subHopf algebra of $A = N' \cap M_1$. 
\end{prop}
\begin{pf}
  Let $S_A$ be an antipode of $A$, i.e., $S_A:A \to A$ is
  an  anti-algebra morphism determined by 
  \begin{displaymath}
    (S_A(a), b) = \overline{(a^*,b^*)} \ \text{for } \
    \forall a \in A \ \text{and} \ \forall b \in 
    B = M' \cap M_2. 
  \end{displaymath}
  Since
  $Be_NB = N' \cap M_2$ by the assumption, 
  for any $a \in A$, there exist $x_i, y_i \in B$ 
  such that $ a = \sum_i x_ie_N y_i$.
  Then $S_A(a) = \sum_i y_ie_N x_i$ 
  (see for example \cite{Szym:Hopf}).
  By the assumption and Lemma \ref{lem:pairing},
  $H$ and $B_{e_{K_1}} = e_{K_1}Be_{K_1}$ 
  are the dual pair of  Hopf algebras with a pairing 
  \begin{displaymath}
    (h , k) = [M:N]^2\tau(h e_Me_N k) \ 
    \text{for } \ \forall h \in H \ 
    \text{and} \ \forall k \in B_{e_{K_1}}.
  \end{displaymath}
  By the fact that 
  $\varphi (e_K) 
  = \frac{[M:N]}{[M:K]}e_{K_1}e_Ne_{K_1}$, 
  for $h \in H $, there exist $s_n, t_n \in B_{e_{K_1}}$ 
  such that $he_{K_1} = \varphi(h) = \sum_n s_n e_N t_n$, 
  where $\varphi$ is defined in Lemma \ref{lem:isom}.
  Then   for $\forall b \in B$, we have
  \begin{align*}
    (S_A(h), b) 
    &= \overline{(h^*, b^*)} \\
    &= [M:N]^2 \tau(b e_N e_M h) \\
    &= [M:N]^2 \sum_n \tau(b e_N e_M s_n e_N t_n) \\
    &= [M:N]^2 \sum_n \tau(b E_{M_1'}^{M'}( e_M s_n) e_N t_n) \\
    &= [M:N]^2 \sum_n \tau(e_Ms_n)\tau(b e_N t_n) \\
    &= [M:N]\sum_n\tau(e_Ms_n)\tau(bt_n).
  \end{align*}
  Since $S_H(h)e_{K_1} = S_{H_{e_{K_1}}}(he_{K_1})
  = \sum_n t_ne_Ns_n$ by the fact that 
  $e_{K_1} \in H'$, 
  we have, for $\forall b \in B$, 
  \begin{align*}
    (S_H(h),b)
    &= [M:N]^2 \tau(S_{H_{e_{K_1}}}(he_{K_1})e_Me_N b)\\
    &= [M:N]^2 \sum_n \tau(t_ne_Ns_ne_Me_N b) \\
    &= [M:N]\sum_n \tau(s_ne_M)\tau(t_nb).
  \end{align*}
  Therefore we have $S_A(h) = S_H(h) \in H$, i.e., 
  $S_A(H) \subset H$.
  
  Let $\Delta_A$ be a comultiplication of $A$, i.e., 
  $\Delta_A:A \to A \otimes A$ is  determined by 
  \begin{displaymath}
    (a, b_1b_2) = (\Delta_A(a), b_1 \otimes b_2) \ 
    \text{for} \ \forall b_1, b_2 \in B. 
  \end{displaymath}
  For $h \in H $, we denote $\Delta_A(h)$  
  by  $\sum_{(h)}h_{(1)}\otimes h_{(2)}$. 
  Since 
  $
    e_M = e_{K_1}e_M 
  $ 
  \ and 
  $
    e_{K_1}h = h e_{K_1}
  $,
  we have 
  \begin{equation}
    \label{center}
    \begin{split}
      (h, b)
      &= [M:N]^2\tau(he_{K_1}e_Me_Nb) \\
      &= [M:N]^2\tau(he_Me_Nbe_{K_1}) \\
      &= (h, b e_{K_1})
      \ \text{for} \ \forall h \in H \
      \text{and} \ \forall b \in B.
    \end{split}
  \end{equation}
  Since $e_{K_1}$ is an element of the center of 
  $B$ by the 
  proof of Theorem \ref{thm:depth2}, we have 
  \begin{align*}
    (h , b_1b_2) 
    &= (h, b_1e_{K_1}b_2e_{K_1})\\
    &= (\Delta_A(h), b_1e_{K_1} \otimes b_2e_{K_1})\\
    &= \sum_{(h)}(h_{(1)}, b_1e_{K_1})(h_{(2)}, b_2e_{K_1})\\
    &= \sum_{(h)}[M:N]^2\tau(e_{K_1}h_{(1)}e_Me_N b_1)
    [M:N]^2\tau(e_{K_1}h_{(2)}e_Me_N b_2)\\
    &= \sum_{(h)}(E_{K_1}^{M_1}(h_{(1)}), b_1)(E_{K_1}^{M_1}(h_{(2)}),b_2), 
    \ \text{for } \ \forall b_1, b_2 \in B.
  \end{align*}
  Since $\sum _{(h)}S_A(h_{(1)})e_Mh_{(2)} \in B$ by
  \cite{Szym:Hopf},
  we have 
  \begin{displaymath}
    \Delta_A(H) \subset H \otimes H
  \end{displaymath}
  by Lemma \ref{lem:A}.
  We have thus proved the theorem.
\end{pf}
\begin{thm}\label{thm:main}
  Let $N \subset M$ be an irreducible, depth 2 inclusion 
  of type II$_1$ factors with finite index and 
  $K$  an intermediate subfactor for $N \subset M$.
  Let $N \subset M \subset M_1 \subset M_2$ and 
  $K \subset M \subset K_1 \subset K_2$ be the Jones 
  towers for $N \subset M$ and $K \subset M$, respectively.
  Then the  depth of $K \subset M$ is 2 if and only if 
  $e_{K_1}$ is an element of the center of $M' \cap M_2$,
  where $e_{K_1}$ is the Jones projection for 
  $K_1 \subset M_1$.
\end{thm}
\begin{pf}
  Suppose that the depth of $K \subset M$ is 2. 
  Then by the proof of Theorem \ref{thm:depth2},
  $e_{K_1}$ is an element of the center of $M' \cap M_2$.
  
  Conversely, suppose that $e_{K_1}$ is an element of  
  the center of $M' \cap M_2$.
  Then for any $h \in H = K' \cap K_1$, we have 
  \begin{align*}
    (S_A(h), b) 
    &= \overline{(h^*, b^*)} \\
    &= [M:N]^2\tau(b^*e_Ne_Mh^*) \\
    &= [M:N]^2\tau(e_{K_1}b^*e_{K_1}e_Ne_Mh^*) \\
    &= (S_A(h), e_{K_1}be_{K_1})
    \quad
    \text{for } \forall b \in B = M' \cap M_2
  \end{align*}
  and hence 
  $
    S_A(H) \subset H
  $.
  Similarly, for any $h \in H$, we have 
  \begin{align*}
    (\Delta_A(h), x \otimes y) 
    &= (h , xy) \\ 
    &= (h, e_{K_1}xe_{K_1}ye_{K_1})\\
    &= (\Delta_A(h), e_{K_1}xe_{K_1} \otimes e_{K_1}ye_{K_1}) 
    \quad
    \text{for }  \forall x, y \in M' \cap M_2,
  \end{align*}
  and hence $\Delta_A(H) \subset H \otimes H$.
  Therefore $H$ is a subHopf algebra of $N' \cap M_1$. 
  By Proposition \ref{prop:subalgebra}, we have 
  $K = M^H$. So the depth of $K \subset M$ is 2.  
\end{pf}
\begin{cor}\label{cor:depth}
  Let $N \subset M$ be an irreducible, depth 2 inclusion 
  of type II$_1$ factors with finite index and 
  $K$ an intermediate subfactor for $N \subset M$.
  Let $N \subset M \subset M_1 \subset M_2$ and 
  $K \subset M \subset K_1 \subset K_2$ be the Jones 
  towers for $N \subset M$ and $K \subset M$, respectively.
  The depth of $N \subset K$ is 2 if and only if 
  $e_K$ is an element of the center of $N' \cap M_1$,
  where $e_K$ is the Jones projection for 
  $K \subset M$.
\end{cor}
\begin{pf}
  Let $K_{-1}$ and $N_{-1}$ be the tunnel constructions  
  for $ N \subset K$ and $N \subset M$, respectively.
  Then the depth of $N_{-1} \subset N$ is 2 and, 
  the depth of $N \subset K$ is 2 if and only if 
  the depth of $K_{-1} \subset N$ is 2.
  And hence, by Theorem \ref{thm:main}, 
  we get the corollary.
\end{pf}
\begin{thm}
  Let $N \subset M$ be an irreducible, depth 2 inclusion 
  of type II$_1$ factors with finite index and 
  $K$  an intermediate subfactor for $N \subset M$.
  Then $K$ is a normal intermediate subfactor of  $N \subset M$ 
  if and only if the depths of $N \subset K$ and $K \subset M$
  are both 2.
\end{thm}
\begin{pf}
  This immediately follows from  
  Theorem \ref{thm:main} and Corollary \ref{cor:depth}.
\end{pf}
\begin{thm}\label{thm:nH}
  Let $N \subset M$ be an irreducible, depth 2 inclusion 
  of type II$_1$ factors with finite index and 
  $K$  an intermediate subfactor for $N \subset M$.
  Then $K$ is a normal intermediate subfactor of 
  $N \subset M$ if and only if 
  $K' \cap K_1$ is a normal subHopf algebra of 
  $N' \cap M_1$, where $K_1$ and $M_1$ are the 
  basic extensions for $N \subset M$ and $K \subset M$,
  respectively.
\end{thm}
\begin{pf}
  Suppose that  $K$ is a normal intermediate subfactor of 
  $N \subset M$. Then $H =  K' \cap K_1$ is 
  a subHopf algebra of 
  $A = N' \cap M_1$ by Proposition \ref{prop:subHopf}.
  Let $\varepsilon_H$ is a counit of $H$. Then
  \begin{displaymath}
    xe_K = \varepsilon_H(x)e_K \quad 
    \text{for }x \in H.
  \end{displaymath}
  Therefore $H^+ = H \cap \ker \varepsilon_H = H (1 - e_K)$.
  Since $(1 - e_K)$ is an element of the center of $A$ 
  by the assumption, we have 
  \begin{displaymath}
    H^+A = AH^+.
  \end{displaymath}  
  Hence $H$ is a normal subHopf algebra of $A$ 
  by Proposition \ref{prop:normalHopf}.
  Conversely, we suppose that 
  $H$ is a normal subHopf algebra of $A$.
  Then by Proposition \ref{prop:subHopf} and 
  Proposition \ref{prop:normalHopf}, 
  $e_K$ and $e_{K_1}$ are elements of the 
  centers of $N' \cap M_1$ and $M' \cap M_2$, 
  respectively and hence $K$ is 
  a normal intermediate subfactor of $N \subset M$.
\end{pf}

\subsection{Lattices of normal intermediate subfactors}
\label{sec:lattice}

Let $N \subset M$ be an irreducible, depth 2  
inclusion of type II$_1$ factors with finite index.
In this subsection we shall prove that 
the set of all normal intermediate subfactors of 
the inclusion $N \subset M$, denoted by 
$\N(N \subset M)$, is a sublattice of 
$\Lat(N \subset M)$.
Moreover, $\N(N \subset M)$ is a modular lattice.
\begin{lem}
  Let $L$ and $K$ be intermediate subfactors of 
  $N \subset M$ and $L_1$ and $K_1$ the basic 
  extensions for $L \subset M$ and $K \subset M$, 
  respectively.  
  Then the basic extension $(L \wedge K)_1$ 
  for $(L \wedge K) \subset M$ is 
  $L_1 \vee K_1$ and
  the basic extension $(L \vee K)_1$ 
  for $(L \vee K) \subset M$ is 
  $L_1 \wedge K_1$.
\end{lem}
\begin{pf}
  By the fact that $(L \cap K)' = (L' \cup K')''$, we have 
  \begin{displaymath}
    (L \wedge K)_1 = J (L \wedge K)' J = L_1 \vee K_1.
  \end{displaymath}
  Similarly, by the fact that 
  $(L \cup K)' = L' \cap K'$, we have 
  \begin{displaymath}
    (L \vee K)_1 = J(L \cup K)' J = L_1 \wedge K_1.
  \end{displaymath}
\end{pf}
We note that if we denote  by $e_A$ 
the Jones projection for $A \subset M$, 
then for $L, K \in \Lat(N \subset M)$, 
we have $e_{L \wedge K} =  e_L \wedge e_K$. 
But $e_{L \vee K} \not= e_L \vee e_K$ in general 
(see \cite{SW:angle}).
\begin{thm}\label{thm:nlattice}
  Let $N \subset M$ be an irreducible, depth 2 
  inclusion of type II$_1$ factors with finite index.
  Then the set of all normal 
  intermediate subfactors $\N(N \subset M)$ is a sublattice of 
  $\Lat(N \subset M)$
\end{thm}
\begin{pf}
  Let $L$ and $K$ be normal intermediate subfactors of 
  $N \subset M$.
  Since $e_L$ and $e_K$ are elements of the center 
  of $N' \cap M_1$ by the assumption, 
  we have $e_{L \wedge K} = e_L \wedge e_K 
  \in \cal Z(N' \cap M_1)$ by the above argument. 
  Observe that 
  \begin{displaymath}
    (L \vee K)' \cap (L \vee K)_1 = 
    (L' \cap L_1) \cap (K' \cap K_1). 
  \end{displaymath}
  Since $L' \cap L_1$ and $K' \cap K_1$ 
  are invariants under the left and right adjoint 
  action of $N' \cap M_1$ 
  (see Definition \ref{defn:defHopf}), 
  so is $(L \vee K)' \cap (L \vee K)_1$.
  Therefore we can see that 
  $(L \vee K)' \cap (L \vee K)_1$ is a 
  normal subHopf algebra $N' \cap M_1$ by 
  the definition.
  Since  $L \vee K$ is a normal intermediate subfactor
  of $N \subset M$ by Theorem \ref{thm:nH}, 
  we have $e_{L \vee K} \in \cal Z(N' \cap M_1)$.
  Applying the same argument for the dual inclusion 
  $M \subset M_1$, we conclude that 
  $L \wedge K$ and $L \vee K$ are normal intermediate 
  subfactors of $N \subset M$.
\end{pf}
\begin{cor}
  Let $N \subset M$ be an irreducible, depth 2 inclusion 
  of type II$_1$ factors.
  Then $\N(N \subset M)$ is a modular lattice.
\end{cor}
\begin{pf}
  This immediately follows from 
  Proposition \ref{prop:quasi}, Theorem \ref{thm:nlattice} 
  and \cite[Theorem 3.9]{Wa:lattice}.
\end{pf}

\begin{thm}
  Let $N \subset M$ be an irreducible, depth 2 inclusion 
  of type II$_1$ factors with finite index.
  Then every maximal chain from $M$ to $N$ in 
  $\N(N \subset M)$ has the same length, i.e., 
  for $A_i (i = 1,2, \dots, m)$, 
  $B_j (j = 1,2, \dots, n) \in \N(N \subset M)$,
  if 
  \begin{displaymath}
    M = A_0 > A_1 > \cdots > A_m = N 
  \end{displaymath}
  and 
  \begin{displaymath}
    M = B_0 > B_1 > \cdots > B_n = N, 
  \end{displaymath}
  then $m = n$, 
  where  
  $X > Y$ 
  means $X \supset Y$ and $X \supseteq K \supseteq Y $, 
  implies $K = X$ or $K = Y$ for $X, Y, K  \in \N(N \subset M)$.
\end{thm}
\begin{pf}
  Since we have the Jordan-Dedekind chain condition holding 
  in modular latteces, this immediately 
  follows from the previous corollary.
\end{pf}
\begin{exmp}\label{exmp:length}
  We denote by  $S_n$ the symmetric group on $n$ letters, 
  $x_1, x_2, \cdots, x_n$ 
  and  $\sigma = (1, 2, 3, \cdots, n)$   
  the  element of 
  $S_n$ with order $n$ and 
  $\<\sigma\>$ the cyclic group generated by $\sigma$.
  Let $\gamma:S_n \to \Aut(P)$ be an outer action of 
  $S_n$ on a type II$_1$ factor $P$ and let 
  \begin{displaymath}
    N = P^{\gamma_{\sigma}} 
    \subset M = P \rtimes_{\gamma}S_{n-1}.
  \end{displaymath}
  Then we can see that 
  $
  S_n = S_{n-1}\<\sigma\> = \<\sigma\>S_{n-1}
  $ 
  and 
  $
  S_{n-1} \cap \<\sigma\> = \{e\}
  $.
  Therefore the depth of $N \subset M$ is 2 
  (see \cite{Sano:com,Yag:vector}).
  We put $K = P \rtimes_{\gamma}A_{n-1}$, 
  where $A_{n-1}$ is the alternating group consists of 
  the even permutations on $x_1, x_2, \dots, x_{n-1}$.
  If $n$ is odd, then the length of 
  $\N(N \subset M)$ is 3 and 
  if $n$ is even, then that is 2
  (we shall show this fact later 
  in Example \ref{exmp:normal}).
\end{exmp}

\section{Some examples}

In this section we shall give some examples of normal 
intermediate subfactors and non normal ones.

\subsection{Group type inclusions}
\label{sec:Majid}
Let $\gamma : G \to \Aut(P)$ be an outer action of a 
discrete group $G$ on a type II$_1$ factor.
Let $A$ and $B$ be finite subgroups of $G$
such that $A \cap B = \{e\}$.
Let $N$ be the fixed point algebra $P^{(A, \gamma)}$ 
and $M$ the crossed product $P \rtimes_{\gamma} B$.
Then $N \subset M$ is an irreducible inclusion
by \cite{BH:comp} and $P$ is normal in $N \subset M$
by Theorem \ref{thm:depth2}.
In this subsection we consider inclusions of this type.

\begin{prop}\label{prop:BH}
  With the above notation, 
  let $H$ be a subgroup of $B$ and $K$ the 
  crossed product $P \rtimes H$.
  Then $K$ is normal in $N \subset M$ if 
  and only if $H$ is a normal subgroup of $B$ 
  and $AH  \cap BA = AH \cap HA$.
\end{prop}
\begin{pf}
  Let $\alpha = \sideset{_N}{_P}{L^2(P)}$ and 
  $\beta = \sideset{_P}{_M}{L^2(M)}$.
  Let $\beta_1 = \sideset{_P}{_K}{L^2(K)}$ and 
  $\beta_2 = \sideset{_K}{_M}{L^2(M)}$.
  Then we have
  \begin{align*}
    \overline{\alpha}\alpha 
    &= \oplus_{a \in A}\sideset{_P}{_P}{L^2(\gamma_a)} \\
    \beta_1\overline{\beta}_1 
    &= \oplus_{h \in H} \sideset{_P}{_P}{L^2(\gamma_h)} \\
    \beta\overline{\beta} 
    &= \oplus_{b \in B} \sideset{_P}{_P}{L^2(\gamma_b)}, 
  \end{align*}
  as in Example \ref{exmp:group}.
  Since $A \cap B = \{e\}$,  
  we have 
  $$
  \left(ab = a'b', \ a, a' \in A, \  
    b, b' \in B
  \right)
  \Longleftrightarrow 
  \left( a = a' \ \text{and } \ b = b'  
  \right).
  $$
  Therefore if 
  $\rho = \sideset{_N}{_M}{L^2(M)} (= \alpha\beta$), 
  then  
  \begin{align*}
    \<\alpha\beta_1\overline{(\alpha\beta_1)}, 
    \rho\overline{\rho} \>
    &= \<\alpha\beta_1\overline{\beta}_1\overline{\alpha}, 
    \alpha\beta\overline{\beta}\overline{\alpha} \> \\
    &= \< \overline{\alpha}\alpha
    \beta_1\overline{\beta}_1, 
    \beta\overline{\beta}\overline{\alpha}\alpha\> \\
    &= \sideset{^{\#}}{}{(AH \cap BA)}
  \end{align*}
  and 
  \begin{align*}
    \<\alpha\beta_1\overline{(\alpha\beta_1)}, 
    \alpha\beta_1\overline{(\alpha\beta_1)} \>
    &= \<\overline{\alpha}\alpha \beta_1\overline{\beta}_1, 
    \beta_1\overline{\beta}_1\overline{\alpha}\alpha\>  \\
    &= \sideset{^{\#}}{}{(AH \cap HA)}.  
  \end{align*}
  Hence $e_K \in \cal Z(N' \cap M_1)$ 
  if and only if 
  $(AH \cap BA) = (AH \cap HA)$ 
  by Proposition \ref{lem:defbi}.
  Suppose $K$ is normal in $N \subset M$. 
  Then $K$ is also normal in $P \subset M$ by Proposition
  \ref{prop:int}.
  Therefore $H$ is a normal subgroup of $B$ 
  by Proposition \ref{prop:crossed}.
  Conversely, if $H$ is a normal subgroup of $B$, 
  i.e., the depth of $K \subset M$ is 2, then we have 
  \begin{align*}
    \<\overline{\beta}_2\beta_2, \overline{\rho}\rho \>
    &= \<\overline{\beta}_2\beta_2, 
    \overline{\beta}_2\overline{\beta}_1\overline{\alpha}
    \alpha\beta_1\beta_2 \> \\
    &= \<\beta_2\overline{\beta}_2\beta_2\overline{\beta}_2, 
    \overline{\beta}_1\overline{\alpha}\alpha\beta_1 \> \\
    &= [B:H] \<\beta_2\overline{\beta}_2, 
    \overline{\beta}_1\overline{\alpha}\alpha\beta_1 \> \\
    &= [B:H]
    = \<\overline{\beta}_2\beta_2,\overline{\beta}_2\beta_2\>.
  \end{align*}
  This proves the proposition.
\end{pf}

Let $G$ be a finite group with two subgroups $A, B$ satisfying 
$G = AB$ and $A \cap B = \{e\}$.
By the uniqueness of the decomposition of an element in 
$G = AB = BA$, 
we can represent $ab$ for $a \in A, b \in B$ 
as
\begin{displaymath}
  ab = \alpha_a(b) \beta_{b^{-1}}(a^{-1})^{-1} \in BA.
\end{displaymath}
Then the matched pair  $(A,B,\alpha, \beta)$ appears 
(see for example \cite{Sano:com}).

\begin{prop}\label{prop:depth2}
  Let $(A, B, \alpha, \beta)$ be the matched pair defined 
  as above and let
  \begin{displaymath}
    M = P \rtimes_{\gamma}B \supset N = P^{(A, \gamma)} = 
    \{x \in P | \gamma_a(x) = x, \forall a \in A \},
  \end{displaymath}
  where $\gamma$ is an outer action of $G$ on II$_1$ factor $P$.
  Then the depth of $N \subset M$ is 2.
\end{prop}
\noindent
See for a proof \cite{Sano:com,Yag:vector}.

%\begin{pf}
%  Let $\rho_1$ and $\rho_2$ be the $N$-$P$ bimodule 
%  $\sideset{_N}{_R}{L^2(P)}$ and $P$-$M$ bimodule 
%  $\sideset{_P}{_M}{L^2(M)}$, respectively.
%  Then we have
%  \begin{displaymath}
%    \overline{\rho}_1\rho_1 
%    \simeq 
%    \oplus_{a \in A} \sideset{_P}{_P}{L^2(\gamma_a)}
%  \end{displaymath}
%  and 
%  \begin{displaymath}
%    \rho_2\overline{\rho}_2 
%    \simeq
%    \oplus_{b \in B} \sideset{_P}{_P}{L^2(\gamma_b)}
%  \end{displaymath}
%  as in Example \ref{exmp:group}.
%  Since $AB = BA$, we have 
%  \begin{align*}
%    \overline{\rho}_1\rho_1 \rho_2\overline{\rho}_2 
%    &= \sum_{a \in A} \sum_{b \in B} 
%    \sideset{_P}{_P}{L^2(\gamma_a)} \otimes_P 
%    \sideset{_P}{_P}{L^2(\gamma_b)} \\
%    &=  \sum_{a \in A} \sum_{b \in B} 
%    \sideset{_P}{_P}{L^2(\gamma_{ab})} \\
%    &=  \sum_{a \in A} \sum_{b \in B} 
%    \sideset{_P}{_P}{L^2(\gamma_{ba})} 
%    = \rho_2\overline{\rho}_2\overline{\rho}_1\rho_1.
%  \end{align*}
%  Therefore if $\rho$ is the $N$-$M$ bimodule 
%  $\sideset{_N}{_M}{L^2(M)} (= \rho_1 \rho_2)$, 
%  we have 
%  \begin{align*}
%    \rho\overline{\rho}\rho 
%    &= \rho_1\rho_2
%    \overline{\rho}_2\overline{\rho}_1
%    \rho_1\rho_2\\
%    &= \rho_1\overline{\rho}_1\rho_1
%    \rho_2\overline{\rho}_2\rho_2 \\
%    &= \underbrace{(\rho_1 \oplus \cdots \oplus \rho_1)}
%    _{[P:N] \  \text{times}} 
%    \underbrace{(\rho_2 \oplus \cdots \oplus \rho_2)}
%    _{[M:P] \  \text{times}} \\
%    &= \underbrace{\rho \oplus \cdots \oplus \rho}
%    _{[M:N] \  \text{times}}.
%  \end{align*}
%  This implies that the depth of $N \subset M$ is 2.
%\end{pf}

\begin{thm}\label{thm:Majid}
  Let $G$ be a finite group with two subgroups $A, B$ satisfying 
  $G = AB$ and $A \cap B = \{e\}$ and 
  $(A, B, \alpha, \beta)$ the associated matched pair.
  Let $\gamma:G \to \Aut(P)$ be an outer action of $G$ on 
  a type II$_1$ factor $P$ and let
  \begin{displaymath}
    M = P \rtimes_{\gamma}B \supset N = P^{(A, \gamma)} = 
    \{x \in P | \gamma_a(x) = x, \forall a \in A \}.
  \end{displaymath}
  If $H$ is a subgroup of $B$ and $K = P \rtimes_{\gamma}H 
  \in \cal L(N \subset M)$, 
  then $K$ is a normal intermediate subfactor for $N \subset M$ 
  if and only if 
  \begin{enumerate}
  \item $H$ is a normal subgroup of $B$,
  \item $\alpha_a(H) = H, \ \forall a \in A$, i.e., 
    $AH = HA$
  \end{enumerate}
  In particular, if $G$ is a semi direct product 
  $B \rtimes A$, then $K$ is a normal in $N \subset M$ 
  if and only if $H$ is a normal subgroup of $G$.
\end{thm}
\begin{pf}
  Since $BA = AB = G$, we have $(AH \cap BA) = AH$. 
  By Proposition \ref{prop:BH}, we have
  $K$ is normal intermediate subfactor in $N \subset M$ 
  if and only if $H$ is a normal subgroup of $B$ and 
  $(AH \cap HA) = AH$, i.e., $AH = HA$ since 
  $\sideset{^{\#}}{}{HA} = \sideset{^{\#}}{}{AH}$.
\end{pf}
\begin{exmp}\label{exmp:normal}
  Let 
  $
  N = P^{\gamma_{\sigma}} 
  \subset M = P \rtimes_{\gamma}S_{n-1}
  $
  be the irreducible inclusion defined as in 
  Example \ref{exmp:length}.
  The depth of $N \subset M$ is 2 by Proposition \ref{prop:depth2}.
  We put $K = P \rtimes_{\gamma}A_{n-1}$.
  If $n$ is odd, then $\sigma$ is an even permutation  
  and we can see that 
  $A_n = A_{n-1}\<\sigma\> = \<\sigma\>A_{n-1}$. 
  Therefore $K$ is normal in $N \subset M$
  by Theorem \ref{thm:Majid}.
  If $n$ is even, then $\sigma$ is an odd permutation. 
  Since the product of an even and odd permutation 
  in either order is odd, and the product of two odd permutation
  is even, 
  $A_{n-1}\<\sigma\>$ is not subgroup of $S_n$ and hence 
  $A_{n-1}\<\sigma\> \not= \<\sigma\>A_{n-1}$.
  Therefore $K$ is not normal in $N \subset M$ 
  by Theorem \ref{thm:Majid}.
  
  Since $S_{n-1}$ is a maximal subgroup of $S_n$, we have 
  if $\<\sigma^k \>S_{n-1} = S_{n-1}\<\sigma^k\>$, then 
  $\<\sigma^k\> = \<\sigma\>$ or $k = 0 (\mod n)$, 
  i.e., 
  there is no normal intermediate subfactor 
  $K$ of $N \subset M$
  such that $N \subsetneqq K \subsetneqq P$ by 
  Theorem \ref{prop:depth2}.
\end{exmp}
\begin{Rem}
  By Example \ref{exmp:normal}, we have completed 
  the proof of Example \ref{exmp:length}.
\end{Rem}

\subsection{Strongly outer actions  and intermediate subfactors.}
\label{sec:strong}

In this subsection we shall study relations  
between strongly 
outer actions  introduced by Choda and Kosaki \cite{CK:strong}
and normal intermediate subfactors.

Let $N \subset M$ be a pair of type II$_1$ factors, and we set 
\begin{displaymath}
 \Aut(M, N ) = \{ \ \theta \in \Aut(M) \ | \ \theta(N) = N \ \}.
\end{displaymath}
Let 
\begin{displaymath}
  N (= M_{-1}) \subset M (= M_0) \subset M_1 \subset 
  M_2 \subset \cdots
\end{displaymath}
be the Jones tower of the pair $N \subset M$, 
and $e_k (\in M_k)$ the Jones projection for the pair 
$M_{k-2} \subset M_{k-1}$. Then 
each automorphism $\theta \in \Aut(M, N)$ is extended to all 
$M_n$ subject to the condition $\theta(e_i) = e_i$.

\begin{Defn}
  An automorphism $\theta \in \Aut(M,N)$ is said to 
  be strongly outer if the following condition is satisfied 
  for all $k \geq -1$:
  \begin{displaymath}
    a \in M_k \ \text{satisfies} \ ax = \theta(x)a 
    \ \text{for all} \ x \in N \Rightarrow a = 0.
  \end{displaymath}
An action $\alpha$ of a group $G$ into $\Aut(M,N)$ is said to 
be strongly outer if 
$\alpha_g$ is strongly outer for all $g \in G$ except for 
the identity $e$.
\end{Defn}
For $\theta \in \Aut(M,N)$, let $\sideset{_N}{_N}{L^2(\theta)}$
be the $N$-$N$ bimodule as in Example \ref{exmp:aut}.
M.~Choda and H.~Kosaki \cite{CK:strong} gave the next 
characterization of strongly outer automorphisms.
\begin{Thm}
  For $\theta \in \Aut(M,N)$, if 
  $\sideset{_N}{_N}{L^2(\theta)}$ does not 
  appear in the irreducible decomposition of 
  $(\rho\overline{\rho})^k$, $k = 1,2, \dots$, 
  then $\theta$ is strongly outer,
  where $\rho$ is the $N$-$M$ bimodule 
  $\sideset{_N}{_M}{L^2(M)}$.
\end{Thm}
%\begin{pf}
%  Suppose that there exists unitary $u \in M_{k-1}$ 
%  such that $uxu^* = \theta(x), \forall x \in N$ for 
%  some $k > -1$.
%  By Lemma \ref{lem:equivalent}, we have 
%  \begin{displaymath}
%    (\rho\overline{\rho})^k 
%    \simeq \sideset{_N}{_N}{L^2(M_{k-1})} 
%    \simeq \sideset{_N}{_{\theta(N)}}{L^2(M_{k-1})}.
%  \end{displaymath}
%  Since $\sideset{_N}{_{\theta(N)}}{L^2(M_{k-1})} 
%  \simeq \sideset{_N}{_N}{L^2(M_{k-1})} \underset{N}{\otimes}
%  \sideset{_N}{_{\theta(N)}}{L^2(N)}$ and 
%  $\sideset{_N}{_N}{L^2(M_{k-1})} 
%  \succ \sideset{_N}{_N}{L^2(N)}$, 
%  we have $(\rho\overline{\rho})^k  
%  \succ \sideset{_N}{_{\theta(N)}}{L^2(N)}$.
%  Therefore $(\rho\overline{\rho})^k 
%  \not\succ \sideset{_N}{_{\theta(N)}}{L^2(N)}, \ 
%  \forall k > -1 
%  \Rightarrow$  $\theta$ is strongly outer.    
%\end{pf}

\begin{lem}\label{lem:irr}
  Let $B \subset A$ be an irreducible  pair of type II$_1$ 
  factors with finite index.
  Let $\gamma:G \to \Aut(A,B)$ be an outer action of a 
  finite group $G$ 
  and
  $\alpha = \sideset{_B}{_A}{L^2(A)}$.
  If $\overline{\alpha}\alpha \not\succ 
  \sideset{_A}{_A}{L^2(\gamma_g)}$ for all 
  $g \in G$ except for the identity $e$, 
  then $B' \cap (A \rtimes_{\gamma}G) = \Bbb C$.
  In particular, 
  if $\gamma$ is strongly outer, then 
  $B \subset A \rtimes_{\gamma}G$ is irreducible.
\end{lem}

\begin{pf}
  Let $\beta = 
  \sideset{_A}{_{A \rtimes_{\gamma}G}}{L^2(A \rtimes_{\gamma}G)}$
  and $\rho = 
  \sideset{_B}
  {_{A \rtimes_{\gamma}G}}{L^2(A \rtimes_{\gamma}G)}$ 
  $(= \alpha\beta)$.
  Then we have 
  \begin{displaymath}
    \beta\overline{\beta} = 
    \sideset{_A}{_A}{L^2(A \rtimes_{\gamma}G)}
    \simeq \underset{g \in G}\oplus 
    \sideset{_A}{_A}{L^2(\gamma_g)}.
  \end{displaymath}
  Therefore 
  \begin{align*}
    \<\rho, \rho\> 
    &= \<\alpha\beta, \alpha\beta \> \\
    &= \< \overline{\alpha}\alpha, \beta\overline{\beta} \> = 1.
  \end{align*}
  This implies that 
  $B' \cap (A \rtimes_{\gamma}G) = \Bbb C$.
\end{pf}

\begin{prop}\label{prop:bi}
  Let $B \subset A$ be an irreducible  pair of type II$_1$ 
  factors with finite index.
  Let $\gamma:G \to \Aut(A,B)$ be an outer action of a 
  finite group $G$ 
  and
  $\alpha = \sideset{_B}{_A}{L^2(A)}$.
  Then $A$ is normal in $B \subset A \rtimes_{\gamma}G$ 
  if and only if 
  $\overline{\alpha}\alpha\overline{\alpha}\alpha
  \not\succ  \sideset{_A}{_A}{L^2(\gamma_g)}$ for all 
  $g \in G$ except for the identity $e$.
\end{prop}

\begin{pf}
  Suppose that   $\overline{\alpha}\alpha\overline{\alpha}\alpha
  \not\succ  \sideset{_A}{_A}{L^2(\gamma_g)}$ for all 
  $g \in G$ except for the identity $e$.
  Let $\beta = 
  \sideset{_A}{_{A \rtimes_{\gamma}G}}{L^2(A \rtimes_{\gamma}G)}$
  and $\rho = 
  \sideset{_B}
  {_{A \rtimes_{\gamma}G}}{L^2(A \rtimes_{\gamma}G)}$ 
  $(= \alpha\beta)$.
  Since $\beta\overline{\beta} 
  \simeq \underset{g \in G}\oplus 
  \sideset{_A}{_A}{L^2(\gamma_g)}$, we have 
  \begin{align*}
    \< \alpha \overline{\alpha}, \rho\overline{\rho} \> 
    &= \<\alpha \overline{\alpha}, 
    \alpha\beta\overline{\beta}\overline{\alpha} \> \\
    &= \< \overline{\alpha}\alpha\overline{\alpha}\alpha, 
    \beta\overline{\beta} \> \\ 
    &= \< \overline{\alpha}\alpha\overline{\alpha}\alpha, 
    \sideset{_A}{_A}{L^2(A)} \>  
    = \< \alpha\overline{\alpha}, \alpha\overline{\alpha}\>.
  \end{align*} 
  Since $\< \beta\overline{\beta}, \overline{\alpha}\alpha \> = 1$ 
  by Lemma \ref{lem:irr}, we have 
  \begin{align*}
    \<\overline{\beta}\beta, \overline{\rho}\rho \> 
    &= \<\overline{\beta}\beta, \overline{\beta}\overline{\alpha}
    \alpha\beta \> \\
    &= \< \beta\overline{\beta}\beta, 
    \overline{\alpha}\alpha\beta \> \\
    &= \sideset{^\# }{}{G} \<\overline{\beta}\beta, 
    \overline{\alpha}\alpha \> \\
    &= \sideset{^\# }{}{G} 
    = \< \overline{\beta}\beta , \overline{\beta}\beta \>.
  \end{align*}
  Therefore $A$ is normal in $B \subset A \rtimes_{\gamma}G$ 
  by Lemma \ref{lem:defbi}.
  
  Conversely, suppose that 
  $\overline{\alpha}\alpha\overline{\alpha}\alpha \succ 
  \sideset{_A}{_A}{L^2(\gamma_g)}$ 
  for some $g (\not= e) \in G$.
  Then we have 
  \begin{align*}
    \< \alpha\overline{\alpha}, \rho\overline{\rho} \> 
    &= \<\overline{\alpha}\alpha\overline{\alpha}\alpha,
    \beta\overline{\beta} \> \\
    &\gneqq \< \overline{\alpha}\alpha\overline{\alpha}\alpha,
    \sideset{_A}{_A}{L^2(A)} \>
    = \< \alpha\overline{\alpha},\alpha\overline{\alpha} \>.
  \end{align*}
  And hence $A$ is not normal in 
  $B \subset A \rtimes_{\gamma}G$.
\end{pf}

\begin{thm}
  Let $B \subset A$ be an irreducible  pair of type II$_1$ 
  factors with finite index.
  If $\gamma:G \to \Aut(M,N)$ is  a strongly 
  outer action of a 
  finite group $G$, then 
  $A$ and $B \rtimes_{\gamma}G$ are normal 
  intermediate subfactors for 
  the inclusion $B \subset A \rtimes_{\gamma}G$.
\end{thm}
\begin{pf}
  This immediately follows from the previous proposition.
\end{pf}

\begin{exmp}
  Let $B \subset A$ be an inclusion of type II$_1$ factors 
  with the principal graph 
  $E_6$, 
  
  
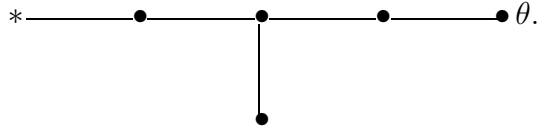
\begin{picture}(200,60)(-115,0)
    \put(8,48){$*$}
    \put(15,50){\line(20,0){40}}
    \put(55,48){$\bullet$}
    \put(61,50){\line(1,0){40}}
    \put(101,48){$\bullet$}
    \put(107,50){\line(1,0){40}}
    \put(147,48){$\bullet$}
    \put(153,50){\line(1,0){40}}
    \put(192,48){$\bullet$}
    \put(103,48){\line(0,-1){35}}
    \put(101,9){$\bullet$}
    \put(200,48){$\theta$.}
  \end{picture}
  
  \noindent
  We put $\alpha = \sideset{_B}{_A}{L^2(A)}$.
  Then we have  
  \begin{displaymath}
    \overline{\alpha}\alpha\overline{\alpha}\alpha
    \succ \sideset{_A}{_A}{L^2(\theta)}.
  \end{displaymath}
  By Proposition \ref{prop:bi},  
  $A$ is not normal in 
  $B \subset A \rtimes_{\theta} \Bbb Z / 2\Bbb Z$.  
\end{exmp}

\ifx\undefined\bysame
\newcommand{\bysame}{\leavevmode\hbox to3em{\hrulefill}\,}
\fi


\begin{thebibliography}{10}

\bibitem{Bi:interm}
D.~Bisch, {\em A note on intermediate subfactors}, Pacific J. Math. {\bf 163}
  (1994), 201--216.

\bibitem{BH:comp}
D.~Bisch and U.~Haagerup, {\em Composition of subfactors: new examples of
  infinite depth subfactors}, preprint.

\bibitem{CK:strong}
M.~Choda and H.~Kosaki, {\em Strongly outer actions for an inclusion of
  factors}, J.\ Funct.\ Anal. {\bf 122} (1994), 315--332.

\bibitem{EK:orbifold}
D.~Evans and Y.~Kawahigashi, {\em Orbifold subfactors from \rm{H}ecke
  algebras}, Comm.\ Math.\ Phys. {\bf 165} (1994), 445--484.

\bibitem{GHJ:coxeter}
F.~Goodman, P.~de~la Harpe, and V.~F.~R. Jones, {\em Coxeter graph and tower of
  algebras}, vol.~14, MSRI Publ, Springer-Verlag, 1989.

\bibitem{Izumi:fusion}
M.~Izumi, {\em Application of fusion rules to classification of subfactors},
  Publ.\ RIMS, Kyoto Univ. {\bf 27} (1991), 953--994.

\bibitem{Izumi:index3}
\bysame, {\em Goldman's type theorem for index 3}, Publ.\ RIMS, Kyoto Univ.
  {\bf 28} (1992), 833--843.

\bibitem{IK:Dn}
M.~Izumi and Y.~Kawahigashi, {\em Classification of subfactors with the
  principal graph \rm{$D_n^{(1)}$}}, J.\ Funct.\ Anal. {\bf 112} (1993),
  257--286.

\bibitem{Jo:index}
V.~F.~R. Jones, {\em Index for subfactors}, Invent.math. {\bf 72} (1983),
  1--25.

\bibitem{Kawahi:flat}
Y.~Kawahigashi, {\em On flatness of \rm{O}cneanu's connection on the dynkin
  diagrams and classification of subfactors}, J.\ Funct.\ Anal. {\bf 127}
  (1995), 63--107.

\bibitem{Ko:charac}
H.~Kosaki, {\em Characterization of crossed product (properly infinite case)},
  Pacific J.\ Math. {\bf 137} (1989), 159--167.

\bibitem{KL:remini}
H.~Kosaki and R.~Longo, {\em A remark on the minimal index of subfactors}, J.\
  Funct.\ Anal. {\bf 107} (1992), 458--470.

\bibitem{Loi:typeIII}
P.~H. Loi, {\em On the theory of index and type \rm{III} factors}, Ph.D.
  thesis, Pennsylvania State Univ., 1988.

\bibitem{Lo:index}
R.~Longo, {\em Index of subfactors and statistics of quantum fields \rm{I}},
  Commun.\ Math.\ Phys. {\bf 126} (1989), 217--247.

\bibitem{Lo:indexII}
\bysame, {\em Index of subfactors and statistics of quantum fields \rm{II}},
  Commun.\ Math.\ Phys. {\bf 130} (1990), 285--390.

\bibitem{mont:Hopf}
S.~Montgomery, {\em Hopf algebras and their actions on rings}, CBMS series
  number 82, 1992.

\bibitem{NT2:Galois2}
M.~Nakamura and Z.~Takeda, {\em On the fundamental theorem of the \rm{Galois}
  theory for finite factors}, Japan Acad. {\bf 36} (1960), 313--318.

\bibitem{NT:Galois}
\bysame, {\em A \rm{Galois} theory for finite factors}, Proc.Japan Acad. {\bf
  36} (1960), 258--260.

\bibitem{Oc:Qg}
A.~Ocneanu, {\em Quantized groups, string algebras and galois theory for
  algebras}, Operator Algebras and Applications, vol.2, London
  Math.~Soc.~Lecture Note Series Vol.136, Cambridge Univ.~Press, 1988,
  pp.~119--172.

\bibitem{Oc:QS}
\bysame, {\em Quantum symmetry, differential geometry of finite graphs, and
  classification of subfactors}, 1991, Univ. of Tokyo Seminary Notes.

\bibitem{PP:entropy}
P.~Pimsner and S.~Popa, {\em Entropy and index for subfactors}, Ann. Sci. Ecole
  Norm. Sup. {\bf 19} (1986), 57--106.

\bibitem{Po:corres}
S.~Popa, {\em Correspondences}, preprint.

\bibitem{Po:reduction}
\bysame, {\em Classification of subfactors: the reduction to commuting
  squares}, Invent.~Math. {\bf 101} (1990), 19--43.

\bibitem{Po:amenable}
\bysame, {\em Classification of amenable subfactors of type \rm{II}}, Acta
  Math. {\bf 172} (1994), 163--255.

\bibitem{Sano:com}
T.~Sano, {\em Commuting co-commuting squares and finite dimensional \rm{Kac}
  algebras}, to appear in Pacific.\ J.\ Math.

\bibitem{SW:angle}
T.~Sano and Y.~Watatani, {\em Angles between two subfactors}, J.\ Operator
  Theory {\bf 32} (1994), 209--242.

\bibitem{Sek:Kac}
Y.~Sekine, {\em Connections associated with finite-dimensional {K}ac algebra
  actions}, Kyushu J. Math. {\bf 49} (1995), 253--269.

\bibitem{Szym:Hopf}
W.~Szyma\'{n}ski, {\em Finite index subfactors and \uppercase{H}opf algebra
  crossecd products}, Proc.\ Amer.\ Math.\ Soc. {\bf 120} (1994), 519--528.

\bibitem{Tamo:normal}
T.~Teruya, {\em A characterization of normal extensions for subfactors}, Proc.\
  Amer.\ Math.\ Soc. {\bf 120} (1994), 781--783.

\bibitem{Wa:lattice}
Y.~Watatani, {\em Lattices of intermediate subfactors}, J. Funct. Anal., to
  appear.

\bibitem{Ya:note}
S.~Yamagami, {\em A note on \rm{Ocneanu's} approach to \rm{Jones'} index
  theory}, Internat.\ J.\ of Math. {\bf 4} (1993), 859--871.

\bibitem{Yag:vector}
\bysame, {\em Vector bundles and bimodules}, Quantum and Non-Commutative
  Analysis, Kluwer Academic, 1993, pp.~321--329.

\bibitem{Yuchi:Kac}
T.~Yamanouchi, {\em Construction of an outer action of a finite-dimensional
  \rm{K}ac algebra on the \rm{AFD} factor of type \rm{II}$_1$}, Inter.\ J.\
  Math. {\bf 4} (1993), 1007--1045.

\end{thebibliography}
\end{document}